\documentclass[a4paper,11pt]{amsart}
\usepackage{amsmath}
\usepackage[T1]{fontenc}
\usepackage{amssymb}
\usepackage{amsthm}
\newcommand{\md}{\,{\rm d}}

\newtheorem{thm}{Theorem}[section]
\newtheorem{prop}[thm]{Proposition}
\newtheorem{lem}[thm]{Lemma}

\newtheorem{proposition}[thm]{Proposition}

\newtheorem{rem}[thm]{Remark}

\numberwithin{equation}{section}
\title[Weyl asymptotics for  magnetic Schr\"odinger operators]
{Weyl asymptotics for  magnetic
Schr\"odinger operators
  and de\,Gennes' boundary condition}

\author{Ayman Kachmar}

\address{A. Kachmar\newline
Universit\'e Paris-Sud\newline D\'epartement de
math\'ematique\newline B\^at.~425\\F-91405 Orsay France}
\email{ayman.kachmar@math.u-psud.fr}
\thanks{This work has been partially supported
by the European Research Network `Postdoctoral Training Program in
Mathematical Analysis of Large Quantum Systems' with contract number
HPRN-CT-2002-00277 and by the ESF Scientific Programme in Spectral
Theory and Partial Differential Equations (SPECT)}
\date{\today}

\begin{document}
\maketitle

\begin{abstract}
This paper is concerned with the discrete spectrum of the
self-adjoint realization of the semi-classical Schr\"odinger operator
with constant magnetic field and associated with the de\,Gennes
(Fourier/Robin) boundary condition. We derive an asymptotic
expansion of the number of eigenvalues below the essential spectrum
(Weyl-type asymptotics).
The methods of proof relies on results concerning the asymptotic
behavior of the first eigenvalue obtained in a previous work
[A.~Kachmar, J.~Math. Phys. {\bf 47} (7) 072106 (2006)].
\end{abstract}

\section{Introduction and main results}

Let $\Omega\subset\mathbb R^2$ be an open domain with regular and
compact boundary. Given a smooth function $\gamma\in
C^\infty(\partial\Omega;\mathbb R)$ and a number
$\alpha\geq\frac12$, we consider the Schr\"odinger operator with
magnetic field~:
\begin{equation}\label{JMP07Op}
P_{h,\Omega}^{\alpha,\gamma} =-(h\nabla-iA)^2,
\end{equation}
whose domain is,
\begin{eqnarray}\label{JMP07DOp}
D\left(P_{h,\Omega}^{\alpha,\gamma}\right)=\big{\{}u\in L^2(\Omega)&:&
(h\nabla-iA)^j\in
L^2(\Omega),~j=1,2,\\
&&\nu\cdot(h\nabla-iA)u+h^\alpha\gamma\,u
=0~{\rm on}~\partial\Omega\big{\}}.\nonumber
\end{eqnarray}
Here $\nu$ is the unit outward normal vector of the  boundary
$\partial\Omega$, $A\in H^1(\Omega;\mathbb R^2)$ is a vector field and
${\rm curl}\,A$ is the magnetic field.
Functions in the domain of $P_{h,\Omega}^{\alpha,\gamma}$ satisfy the
{\it de\,Gennes} boundary condition.\\
 \indent The operator $P_{h,\Omega}^{\alpha,\gamma}$ arises from the analysis
 of the onset of superconductivity for a superconductor placed
 adjacent to another materials. For the physical  motivation and
 the mathematical
 justification of considering this type of boundary condition and not
 the usual Neumann condition ($\gamma\equiv0$), we invite the
 interested reader to see the book of de\,Gennes \cite{deGe} and the
 papers \cite{Kach1, Kach2, Kach3, Kach4}. We would like to mention that when
$\gamma\equiv0$, the operator $P_{h,\Omega}^{\alpha,\gamma}$ has been the
 subject of many papers, see \cite{FoHe} and the references therein.\\
We shall restrict ourselves with the case of {\it constant} magnetic
field, namely when
\begin{equation}\label{JMP07MF=1}
{\rm curl}\,A=1\quad{\rm in }~\overline\Omega.
\end{equation}
It follows from the well-known inequality
\begin{equation}\label{JMP07-Simon}
\int_\Omega|(h\nabla-iA)u|^2\md x\geq h\int_\Omega|u|^2\md
x,\quad\forall~u\in C_0^\infty(\Omega),\end{equation}
and from a `magnetic' Persson's Lemma (cf. \cite{Pe, Bon}), that the
bottom of the essential spectrum of $P_{h,\Omega}^{\alpha,\gamma}$ is
above $h$. Assuming that the boundary of $\Omega$ is smooth and
compact, then it follows from \cite{Kach1} that (in the parameter regime
 $\alpha\geq\frac12$), the operator $P_{h,\Omega}^{\alpha,\gamma}$ has discrete
spectrum below $h$.
Thus, given $b_0<1$, one is led to estimate the size of the discrete
spectrum below $b_0h$, i.e. we look for the asymptotic behavior of
the number
\begin{equation}\label{JMP07N}
N(\alpha,\gamma\,;\,b_0h)
\end{equation}
of eigenvalues of $P_{h,\Omega}^{\alpha,\gamma}$ (taking
multiplicities into account) included in the interval $]0,b_0h]$.\\
For the case with {\it non-constant} magnetic field and {\it Neumann}
boundary condition, this problem has been analyzed by R.~Frank
\cite{Fr} (related questions are also treated in \cite{CV, Iv, Ta, Tr}).
As we shall see, depending on the type of the boundary condition, one
can produce much additional  eigenvalues below the essential spectrum.\\
To state the results concerning $N(\alpha,\gamma\,;\,b_0h)$, we need
to introduce some notation. Let us introduce the {\it smooth}
functions, which arise
from the analysis of the model-operator in the half-plane (see
\cite[Section~II]{Kach1}),
\begin{equation}\label{JMP07Theta}
\mathbb R\times\mathbb R_+\ni(\gamma,\xi)\mapsto\mu_1(\gamma,\xi),\quad
 \mathbb R\ni \gamma\mapsto\Theta(\gamma),
\end{equation}
where
$$\mu_1(\gamma,\xi)=\inf_{u\in B^1(\mathbb R_+),\,u\not\equiv0}
\frac{\displaystyle\int_{\mathbb
    R_+}\left(|u'(t)|^2+|(t-\xi)u(t)|^2\right)\md t+\gamma\,|u(0)|^2}
{\displaystyle\int_{\mathbb R_+}|u(t)|^2\md t},
$$
$$
\Theta(\gamma)=\inf_{\xi\in\mathbb R}\mu_1(\gamma,\xi),$$
and the space $B^1(\mathbb R_+)$ consists of functions in the space
$H^1(\mathbb R_+)\cap L^2(\mathbb R_+; t^2\md t)$.\\
Actually, $\mu_1(\gamma.\xi)$ is the first eigenvalue of the
self-adjoint operator
$$-\partial_t^2+(t-\xi)^2\quad {\rm in}\quad L^2(\mathbb R_+)$$
associated with the boundary condition $u'(0)=\gamma\,u(0)$. The
eigenvalues of this operator form an increasing sequence which we
denote  $\left(\mu_j(\gamma,\xi)\right)_{j\in\mathbb N}$, see
Subsection~\ref{Sub-Sec-M0-1} for more details.\\
 When
$\gamma=0$, we write as in the usual case (see \cite{FoHe})
$$\mu_1(\xi):=\mu_1(0,\xi),\quad\Theta_0:=\Theta(0).$$
Furthermore, we denote by,
\begin{equation}\label{JMP07C1}
C_1(\gamma)=\frac13\left(1+\gamma\sqrt{\Theta(\gamma)+\gamma^2}\right)^2
\Theta'(\gamma),\quad C_1=C_1(0).
\end{equation}

We are ready now to state our main results.

\begin{thm}\label{JMP07thm1}
Assume that  $\alpha>\frac12$ and
$\Theta_0<b_0<1$. Then as $h\to0$,
\begin{equation}\label{JMP07thm11}
N(\alpha,\gamma\,;\,b_0h)=\left(\frac{|\partial\Omega|}{2\pi\sqrt{h}}
\,\big{|}
\{\xi\in\mathbb R:~\mu_1(\xi)<b_0\}\big{|}
\right)
(1+o(1)).
\end{equation}
On the other hand, if $\alpha=\frac12$ and $\Theta(\gamma_0)<b_0<1$,
then as $h\to0$,
\begin{eqnarray}\label{JMP07thm12}
&&N(\alpha,\gamma\,;\,b_0h)=\\
&&\hskip1cm\left(\frac1{2\pi\sqrt{h}} \int_{\partial\Omega}
\sum_{j=1}^\infty\big{|}\{\xi\in\mathbb R:~
\mu_j(\gamma(s),\xi)<b_0\}\big{|} \md
s\right)(1+o(1)).\nonumber\end{eqnarray} Here
\begin{equation}\label{JMP07gamma0}
\gamma_0=\min_{s\in\partial\Omega}\gamma(s).\end{equation}
If we suppose furthermore that $\gamma_0\geq0$,
then (\ref{JMP07thm12}) simplifies to
\begin{eqnarray}\label{JMp07thm12'}
&&N(\alpha,\gamma\,;\,b_0h)=\\
&&\hskip1cm\left(\frac1{2\pi\sqrt{h}} \int_{\partial\Omega}
\big{|}\{\xi\in\mathbb R:~
\mu_1(\gamma(s),\xi)<b_0\}\big{|} \md
s\right)(1+o(1)).\nonumber\end{eqnarray}
\end{thm}

By taking $\gamma\equiv0$, we recover in Theorem~\ref{JMP07thm1} the
result of R.~Frank \cite{Fr}. We notice that when $\alpha=\frac12$
and $\gamma$ is constant, we have additional eigenvalues than the
usual case of Neumann boundary condition if $\gamma<0$ and less
eigenvalues if $\gamma>0$. This is natural as we apply  the
variational min-max  principle. However, when $\gamma_0=0$ or
$\alpha>\frac12$, Theorem~\ref{JMP07thm1} fails to give a comparison
with the Neumann case, i.e. we   have no more information about the
size of the difference:
$$N(\alpha,\gamma\,;\, b_0h)-N(0;\,b_0h).$$
This is at least a
motivation
for  some of the  next results, where we take $b_0=b_0(h)$ asymptotically close
to $\Theta_0$, each time with an appropriate scale (this will cover
also the case $b_0=\Theta_0$).

\begin{thm}\label{JMP07thm2}
If $\frac12<\alpha<1$ 
then for all $a\in\mathbb R$,
\begin{eqnarray}\label{JMP07thm21}
&&N\left(\alpha,\gamma\,;\,h\Theta_0+3aC_1\,h^{\alpha+\frac12}\right)\\
&&\hskip2cm
=\frac1{\pi\sqrt{h^{\frac32-\alpha}\sqrt{\Theta_0}}}\left(\int_{\partial\Omega}
\sqrt{(a-\gamma(s))_+}\,\md s\right)(1+o(1)).\nonumber
\end{eqnarray}
\end{thm}

In the particular case when  
the function $\gamma$ is constant,  the leading order term in
(\ref{JMP07thm21}) will vanish when $a$ is taken equal to $\gamma$.
In this specific regime, Theorem~\ref{JMP07thm4} (more precisely the
formula in (\ref{JMP07thm42})) will substitute
Theorem~\ref{JMP07thm2}.

\begin{thm}\label{JMP07thm3}
Assume that $\alpha=1/2$. Let $0<\varrho<\frac12$,
$\zeta_0>0$, $h_0>0$ and $]0,h_0]\ni h\mapsto c_0(h)\in\mathbb R_+$
a function such
that $\lim_{h\to0}c_0(h)=\infty$. If
$$c_0(h)\,h^{1/2}\leq\left|\lambda-\Theta(\gamma_0)\right|\leq
\zeta_0h^\varrho,\quad\forall~ h\in]0,h_0],$$ then we have the asymptotic
formula,
\begin{eqnarray}\label{JMP07thm31}
 &&N\left(\alpha,\gamma\,;\,h\lambda\right)\\
&&\hskip1cm=
\left(\frac1\pi\int_{\partial\Omega}
\sqrt{\frac{\left[\lambda-\Theta(\gamma(s))\right]_+}
{h\Theta'(\gamma(s))\sqrt{\Theta(\gamma(s))+\gamma(s)^2}}} \,\md
s\right)(1+o(1)),\nonumber
\end{eqnarray}
where the function $\Theta(\cdot)$ being
 introduced in (\ref{JMP07Theta}).
\end{thm}

\begin{rem}\label{generic-gamma}
{\rm Theorem~\ref{JMP07thm3} becomes of particular interest when the
function $\gamma$ has a unique non-degenerate minimum and
$\lambda=\Theta(\gamma_0)+ah^\beta$, for some $a\in\mathbb R_+$ and
$\beta\in]0,\frac12[$. In this case, we have in the support of
$[\lambda-\Theta(\gamma(s))]_+$,
$$\Theta(\gamma(s))=\Theta(\gamma_0)+c_1\, s^2+\mathcal
O(h^{3\beta/2}),$$ for an explicit constant $c_1>0$ determined by
the
functions $\gamma$ and $\Theta$.\\
Therefore, the asymptotic expansion (\ref{JMP07thm31}) reads in this
case, for some explicit constant $c_2>0$,}
\begin{equation}\label{corol-thm3}
N\left(\alpha,\gamma\,;\,h\lambda\right)=c_2\,a\sqrt{a}\,h^{\beta-\frac12}(1+o(1)).
\end{equation}
\end{rem}


The next theorem deals with the regime where the scalar curvature
becomes effective in the asymptotic expansions.

\begin{thm}\label{JMP07thm4}\
\begin{enumerate}
\item Assume that $\alpha=1$.  Then, for all $a\in\mathbb R$, the following
asymptotic expansion holds as $h\to0$,
\begin{eqnarray}\label{JMP07thm41}
&&N\left(1,\gamma\,;\,h\Theta_0+aC_1h^{3/2}\right)=\\
&&\hskip1cm\frac1{\pi\sqrt{3h^{1/2}\sqrt{\Theta_0}}}\left(
\int_{\partial\Omega}\sqrt{(\kappa_{\rm r}(s)-3\gamma(s)+a)_+}\md
  s\right)
(1+o(1)),\nonumber
\end{eqnarray}
where $\kappa_{\rm r}$ is the scalar curvature of $\partial\Omega$.
\item If the function $\gamma$ is constant, then for all
$\alpha>1/2$ and $a\in\mathbb R$, we have the asymptotic expansion,
\begin{eqnarray}\label{JMP07thm42}
&&N\left(\alpha,\gamma\,;\,h\Theta(h^{\alpha-1/2}\gamma)
+aC_1h^{3/2}\right)=\\
&&\hskip2cm\frac1{\pi\sqrt{3h^{1/2}\sqrt{\Theta_0}}}\left(
\int_{\partial\Omega}\sqrt{(\kappa_{\rm r}(s)+a)_+}\md
s\right)(1+o(1)).\nonumber
\end{eqnarray}
\item If the function $\gamma$ is constant and $\alpha=1/2$, then
for all $a\in\mathbb R$,  we have
\begin{eqnarray}\label{JMP07thm43}
&&N\left(\alpha,\gamma\,;\,h\Theta(\gamma)
+aC_1(\gamma) h^{3/2}\right)=\\
&&\hskip1cm\frac{1+\gamma\sqrt{\Theta(\gamma)+\gamma^2}}
{\pi\sqrt{3h^{1/2}\sqrt{\Theta(\gamma)+\gamma^2}}}\left(
\int_{\partial\Omega}\sqrt{(\kappa_{\rm r}(s)+a)_+}\md
s\right)(1+o(1)).\nonumber
\end{eqnarray}
Here $C_1(\gamma)$ has been defined in (\ref{JMP07C1}).
\end{enumerate}
\end{thm}

The proof of Theorems~\ref{JMP07thm1}-\ref{JMP07thm4} is through
careful estimates in the semi-classical regime of the quadratic form
$$u\mapsto q_{h,\Omega}^{\alpha,\gamma}(u)=\int_\Omega|(h\nabla-iA)u|^2\md x+
h^{1+\alpha}\int_{\partial\Omega}\gamma(s)|u(s)|^2\md s\,.$$ These
estimates are essentially obtained in \cite{HeMo} when
$\gamma\equiv0$, then adapted to situations involving the de\,Gennes
boundary condition in \cite{Kach1,KachTh}. We shall follow closely
the arguments of \cite{Fr} but we also require to use various
properties of the function $\gamma\mapsto\Theta(\gamma)$
established in \cite{Kach1}.\\

The paper is organized in the following way. Section~\ref{Sec-MO}
is devoted to the analysis of the model operator in a half-cylinder
when the function $\gamma$ is constant. Section~\ref{proofThm1}
extends the result obtained for the model case in a half-cylinder for
a general domain by which we  prove
Theorem~\ref{JMP07thm1}. Section~\ref{CurvEffects} deals with model
operators
on weighted $L^2$ spaces which serve in proving
Theorems~\ref{JMP07thm2}-\ref{JMP07thm4}.

\section{Analysis of the model operator}\label{Sec-MO}
\subsection{A family of one-dimensional differential operators}
\label{Sub-Sec-M0-1}
Let us recall the main results obtained in \cite{Kach0, Kach1}
concerning a
family of differential operators with Robin boundary condition. Given $(\gamma,\xi)\in\mathbb R\times\mathbb R$, we
define the quadratic form,
\begin{equation}\label{k2-qfk1}
B^1(\mathbb R_+)\ni u\mapsto q[\gamma,\xi](u)=\int_{\mathbb
R_+}\left(|u'(t)|^2+|(t-\xi)u(t)|^2\right)dt+\gamma|u(0)|^2,
\end{equation}
where, for a positive integer $k\in\mathbb N$ and a given interval
$I\subseteq\mathbb R$, the space $B^k(I)$ is defined by~:
\begin{equation}\label{Bk-sp}
B^k(I)=\{u\in H^k(I);\quad t^ju(t)\in L^2(I),\quad \forall
j=1,\cdots, k\}.
\end{equation}
By Friedrichs Theorem, we can associate to the quadratic form
(\ref{k2-qfk1}) a self adjoint operator  $\mathcal L[\gamma,\xi]$
with domain,
$$D(\mathcal L[\gamma,\xi])=\{u\in B^2(\mathbb R_+);\quad
u'(0)=\gamma u(0)\},$$ and associated to the differential operator,
\begin{equation}\label{L-dop}
\mathcal L[\gamma,\xi]=-\partial_t^2+(t-\xi)^2.
\end{equation} We denote by
$\{\mu_j(\gamma,\xi)\}_{j=1}^{+\infty}$ the increasing sequence
of eigenvalues of $\mathcal L[\gamma,\xi]$. When $\gamma=0$ we
write,
\begin{equation}\label{l-ga=0}
\mu_j(\xi):=\mu_j(0,\xi),\quad \forall j\in\mathbb N,\quad
\mathcal L^N[\xi]:=\mathcal L[0,\xi].
\end{equation}
We also denote by $\{\mu^D_j(\xi)\}_{j=1}^{+\infty}$ the
increasing sequence of eigenvalues of the Dirichlet realization of
$-\partial_t^2+(t-\xi)^2$.\\
By the min-max principle, we have,
\begin{equation}\label{l1-ga,xi}
\mu_1(\gamma,\xi)=\inf_{u\in B^1(\mathbb
R_+),u\not=0}\frac{q[\gamma,\xi](u)}{\|u\|^2_{L^2(\mathbb R_+)}}.
\end{equation}
Let us denote by $\varphi_{\gamma,\xi}$ the positive (and
$L^2$-normalized) first eigenfunction of $\mathcal L[\gamma,\xi]$.
It is proved in \cite{Kach1} that the functions
$$(\gamma,\xi)\mapsto\mu_1(\gamma,\xi),\quad (\gamma,\xi)\mapsto
\varphi_{\gamma,\xi}\in L^2(\mathbb R_+)$$ are regular (i.e. of
class $C^\infty$), and we have the following formulae,
\begin{eqnarray}
&&\partial_\xi\mu_1(\gamma,\xi)=-\left(\mu_1(\gamma,\xi)-\xi^2+\gamma^2\right)
|\varphi_{\gamma,\xi}(0)|^2,\label{Th7.5.5}\\
&&\partial_\gamma\mu_1(\gamma,\xi)=|\varphi_{\gamma,\xi}(0)|^2.\label{Th7.5.6}
\end{eqnarray}
Notice that (\ref{Th7.5.6}) will yield that the function
$$(\gamma,\xi)\mapsto\varphi_{\gamma,\xi}(0)$$
is also regular of class $C^\infty$.\\
We define the function~:
\begin{equation}\label{Th-gam}
\Theta(\gamma)=\inf_{\xi\in\mathbb R}\mu_1(\gamma,\xi).
\end{equation}
It is a result of~\cite{DaHe} that there exists a unique
$\xi(\gamma)>0$ such that,
\begin{equation}\label{xi-gam}
\Theta(\gamma)=\mu_1(\gamma,\xi(\gamma)),\quad \Theta(\gamma)<1\,,
\end{equation}
and $\xi(\gamma)$ satisfies (cf. \cite{Kach1}),
\begin{equation}\label{xi(gam)}
\xi(\gamma)^2=\Theta(\gamma)+\gamma^2.
\end{equation}
Moreover, the function $\Theta(\gamma)$ is of class $C^\infty$ and
satisfies,
\begin{equation}\label{Th'-gam}
\Theta'(\gamma)=|\varphi_\gamma(0)|^2,
\end{equation}
where $\varphi_\gamma$ is the positive (and $L^2$-normalized)
eigenfunction associated to $\Theta(\gamma)$~:
\begin{equation}\label{varphi-gam}
\varphi_\gamma=\varphi_{\gamma,\xi(\gamma)}. \end{equation} When
$\gamma=0$, we write,
\begin{equation}\label{Th-0}
\Theta_0:=\Theta(0),\quad \xi_0:=\xi(0).
\end{equation}
It is a consequence of (\ref{Th'-gam}) that the constant  $C_1$
introduced in (\ref{JMP07C1}) can be defined by the alternative manner,
\begin{equation}\label{C1}
C_1:=\frac{|\varphi_0(0)|^2}3.
\end{equation}
Let us recall an important consequence of standard Sturm-Liouville
theory
(cf. \cite[Lemma~2.1]{Fr}).
\begin{lem}\label{RuFr-Lemm-kach}
For all $\xi\in\mathbb R$, we have
$$\mu_2(\xi)>\mu^D_1(\xi)>1.$$
\end{lem}
Let us also introduce,
\begin{equation}\label{Theta2}
\Theta_k(\gamma)=\inf_{\xi\in\mathbb
  R}\mu_k(\gamma,\xi)\,,\quad\forall~k\in\mathbb N.\end{equation}
Another consequence of Sturm-Liouville theory  that we shall need is
the following result on $\Theta_2(\gamma)$.
\begin{lem}\label{HeRu}
For any $\gamma\in\mathbb R$, we have,
$$\Theta_2(\gamma)>\Theta(\gamma).$$
\end{lem}
\paragraph{\bf Proof.}
Let us introduce the continuous function
$f(\gamma)=\Theta_2(\gamma)-\Theta(\gamma)$. Using the min-max
principle, it follows from (\ref{xi-gam})  and
Lemma~\ref{RuFr-Lemm-kach} that $f(0)>0$.
It is then sufficient to prove that the function $f$ never vanish.
Suppose by contradiction that there is some $\gamma_0\not=0$ such
that $\Theta_2(\gamma_0)=\Theta(\gamma_0)$. By the same method used
in \cite{Kach1} one is able to prove that there exists
$\xi_2(\gamma_0)>0$ such that
$$\Theta_2(\gamma_0)=\mu_2(\gamma_0,\xi_2(\gamma_0)),$$
and that $\xi_2(\gamma_0)^2=\Theta_2(\gamma_0)+\gamma_0^2$. Therefore
we get  $\xi_2(\gamma_0)=\xi(\gamma_0)$.\\
Now, by Sturm-Liouville theory, the eigenvalues of the operator
$\mathcal L[\gamma_0,\xi(\gamma_0)]$ are all simple, whereas, by the
above, we get a degenerate eigenvalue
$$\mu_1(\gamma_0,\xi(\gamma_0))=\Theta(\gamma_0)=\mu_2(\gamma_0,\xi(\gamma_0))\,,$$
which is the desired contradiction.\hfill$\Box$\\

One more useful result in Sturm-Liouville theory is the following.

\begin{lem}\label{St-Li-muk}
Let $\gamma\in\mathbb R_-$ and $k\in\mathbb N$. Then
$\Theta_k(\gamma)<2k+1$ and for all $b_0\in]\Theta_k(\gamma),2k+1[$,
the equation
$$\mu_k(\gamma,\xi)=b_0$$
has exactly two solutions $\xi_{k,-}(\gamma,b_0)$ and
$\xi_{k,+}(\gamma,b_0)$. Moreover,
$$\{\xi\in\mathbb
R~:~\mu_k(\gamma,\xi)<b_0\}=\left]-\xi_{k,-}(\gamma,b_0),\xi_{k,+}(\gamma,b_0)
\right[\,.$$
\end{lem}
\paragraph{\bf Proof.}
We can study the variations of the function
$\xi\mapsto\mu_k(\gamma,\xi)$ using exactly the same method of
\cite{Kach1, KachTh, DaHe}. We obtain that the function
$\xi\mapsto\mu_k(\gamma,\xi)$ attains a unique non-degenerate
minimum at the point
$\xi_k(\gamma)=\sqrt{\Theta_k(\gamma)+\gamma^2}$\,, and  analogous
formulae to (\ref{Th7.5.5})-(\ref{Th7.5.6}) continue to hold for
$(\gamma,\xi)\mapsto\mu_k(\gamma,\xi)$. Moreover,
$\displaystyle\lim_{\xi\to-\infty}\mu_k(\gamma,\xi)=\infty$ and
$\displaystyle\lim_{\xi\to
\infty}\mu_k(\gamma,\xi)=2k+1$.\\
For instance, the restrictions of the function
$\xi\mapsto\mu_k(\gamma,\xi)$ to the intervals
$]-\infty,\xi_k(\gamma)[$
    and $]\xi_k(\gamma),\infty[$ are invertible.\hfill$\Box$\\

It is a result of the variational min-max principle that the
function $\gamma\mapsto\Theta_k(\gamma)$ is continuous, see
\cite[Proposition~2.5]{Kach1} for the case $k=1$. Thus the set
\begin{equation}\label{eq-Uk}
U_k=\{(\gamma,b)\in\mathbb R\times\mathbb
R~:~\Theta_k(\gamma)<b<2k+1\}\quad {\rm is~open~in~}\mathbb
R^2.\end{equation}

\begin{lem}\label{reg-xi-k}
The functions
$$U_k\ni(\gamma,b)\mapsto \xi_{k,\pm}(\gamma,b)$$
admit continuous extensions
$$\mathbb
R\times]-\infty,2k+1[\,\mapsto\overline\xi_{k,\pm}(\gamma,b)\,.$$
\end{lem}
\paragraph{\bf Proof.}
Using the regularity of $\mu_k(\gamma,\xi)$, the implicit function
theorem applied to
$$U_k\times\mathbb R\ni(\gamma,b,\xi)\mapsto\mu_k(\gamma,\xi)-b$$
near  $(\gamma_0,b_0,\xi_{k,\pm}(\gamma_0,b_0))$ (for an arbitrary
point $(\gamma_0,b_0)\in U_k$) permits to deduce that the functions
$$U_k\ni(\gamma,b)\mapsto \xi_{k,\pm}(\gamma,b)\quad{\rm
are~}C^1\,.$$ We then define the following continuous extensions of
$\xi_{k,\pm}$,
$$\overline\xi_{k,\pm}(\gamma,b)=\left\{
\begin{array}{l}
\xi_{k,\pm}(\gamma,b)\,,\quad{\rm if~}\Theta_k(\gamma)<b<2k+1\,,\\
\xi_k(\gamma)\,,\quad{\rm if~}\Theta_k(\gamma)\geq b\,,
\end{array}\right.
$$
where $\xi_k(\gamma)$ is the unique non-degenerate minimum of
$\xi\mapsto\mu_k(\gamma,\xi)$.\hfill$\Box$\\

The next lemma justifies that the sum on the right hand side of
(\ref{JMP07thm12}) is indeed finite.

\begin{lem}\label{justification}
For each $M>0$ and $b_0\in]0,1[$, there exists a constant $C>0$ such
that, for all $\gamma\in]-M,M[$ and $b\in]\Theta(\gamma),b_0[$, we
have
$$\sum_{j=1}^\infty\big{|}\{\xi\in\mathbb R~:~\mu_j(\gamma,\xi)<b\}\big{|}
\leq C.$$
\end{lem}
\paragraph{\bf Proof.}
Let us notice that for all $j\geq1$,
$$\{\xi\in\mathbb R~:~\mu_j(\gamma,\xi)<b\}\subset
\{\xi\in\mathbb R~:~\mu_1(\gamma,\xi)<b_0\}\,,$$ and for all
$\gamma\in]-M,M[$ (using the monotonicity of
$\eta\mapsto\mu_1(\eta,\xi)$),
$$\{\xi\in\mathbb
R~:~\mu_1(\gamma,\xi)<b_0\}\subset\{\xi\in\mathbb
R~:~\mu_1(-M,\xi)<b_0\}\,.$$ Consequently, there exists a constant
$\widetilde M>0$  such that
$$\{\xi\in\mathbb R~:~\mu_j(\gamma,\xi)<b\}\subset [-\widetilde M,\widetilde
  M],\quad\forall~b\leq b_0\,,~
\forall~\gamma\in]-M,M[.$$
Since the functions
$$\xi\mapsto\mu_j(\gamma,\xi)\quad (j\in\mathbb N)$$
are regular, we introduce constants
$\left(\xi_j(M)\right)_{j\in\mathbb N}\subset[-\widetilde
M,\widetilde M]$ by
$$\mu_j(-M,\xi_j(M))=
\min_{\xi\in[-\widetilde M,\widetilde M]}\mu_j(-M,\xi)\,.$$ We claim
that
\begin{equation}\label{claim-j}
\lim_{j\to\infty}\mu_j(-M,\xi_j(M))=\infty.\end{equation} Once this
claim is proved, we get the result of the lemma, since by
monotonicity
$$\mu_j(\gamma,\xi)\geq \mu_1(-M,\xi)\quad\forall~\gamma\geq-M\,,~\forall~\xi\in\mathbb R .$$
Let us assume by contradiction that the claim (\ref{claim-j}) were
false. Then we may find a constant $\mathcal M>0$ and a subsequence
$(j_n)$  such that
\begin{equation}\label{claim-j*}
\mu_{j_n}(-M,\xi_{j_n}(M))\leq\mathcal M\,,\quad\forall~n\in\mathbb N\,.
\end{equation}
Since $-\widetilde M\leq \xi_{j_n}(M)\leq \widetilde M$ for all $n$, we
get a subsequence, denoted again by $\xi_{j_n}(M)$, such that
$$\lim_{n\to\infty}\xi_{j_n}(M)=\zeta(M)\in[-\widetilde
  M,\widetilde M].
$$
It is quiet easy, by comparing the corresponding quadratic forms, to
prove the existence of a constant $C>0$ such that, for all
$\varepsilon\in]0,\frac12[$ and $n\in\mathbb N$, we have the estimate
\begin{eqnarray}\label{mu-j->infty}
\mu_{j_n}(-M,\xi_{j_n}(M))&\geq&(1-\varepsilon)\mu_{j_n}(-2M,\zeta(M))\\
&&-C\left(\varepsilon+\varepsilon^{-1}|\xi_{j_n}(M)-\zeta(M)|^2\right).\nonumber
\end{eqnarray}
We shall provide some details concerning the above estimate, but we
would like first to achieve the proof of the lemma. Notice that,
since the operator $\mathcal L[-M,\zeta(M)]$ has compact resolvent,
then
$$\lim_{n\to\infty}\mu_{j_n}(-2M,\zeta(M))=\infty.$$
Upon choosing $\varepsilon=|\xi_{j_n}(M)-\zeta(M)|$, we get from
(\ref{mu-j->infty}) that
$$\lim_{n\to\infty}\mu_{j_n}(-M,\xi_{j_n}(M))=\infty\,,$$
contradicting thus (\ref{claim-j*}).\\
We conclude by some wards concerning the proof of
(\ref{mu-j->infty}). Notice that, for a normalized $L^2$-function
$u$, we have by Cauchy-Schwarz inequality:
\begin{eqnarray*}
2\left|\int_0^\infty(\zeta-\xi_{j_n}(M))(t-\zeta)|u|^2\md
t\right|&\leq&2|\zeta-\xi_{j_n}(M)|\times\|(t-\zeta)u\|_{L^2(\mathbb
R_+)}\\
&\leq&\varepsilon\|(t-\zeta)u\|^2_{L^2(\mathbb
R_+)}+\varepsilon^{-1}|\zeta-\xi_{j_n}(M)|^2\,, \end{eqnarray*} for any
$\varepsilon>0$. On the other hand, writing
$$(t-\xi_{j_n}(M))^2=(t-\zeta)^2+(\xi-\xi_{j_n}(M))^2+2(\zeta-\xi_{j_n}(M))
(t-\zeta)\,,$$
we get the following comparison of the quadratic forms
$$q[-M,\xi_{j_n}(M)](u)\geq
q[-M,\zeta](u)-\varepsilon\|(t-\zeta)u\|_{L^2(\mathbb
R_+)}^2-\varepsilon^{-1}|\zeta-\xi_{j_n}(M)|^2\,,$$ where
$q[-M,\cdot]$ has been introduced in (\ref{k2-qfk1}). Noticing that
for $\varepsilon\in]0,\frac12[$, $\frac{-M}{1-\varepsilon}\geq
  -2M$, the
application of the min-max principle permits then to conclude the
desired bound (\ref{mu-j->infty}).
\hfill$\Box$\\

\begin{rem}\label{rem-justification}
Once the asymptotic expansion (\ref{JMP07thm12}) is proved, the formula
(\ref{JMp07thm12'}) becomes a consequence of
Lemma~\ref{RuFr-Lemm-kach}.
\end{rem}

The next lemma will play a crucial role in establishing the main results
of this paper.

\begin{lem}\label{lem-referee}
The function
$$\mathcal S:\mathbb R\times]-\infty,1[\,\ni(\gamma,b)\mapsto \sum_{j=1}^\infty\left|
\{\xi\in\mathbb R~:~\mu_j(\gamma,\xi)<b\}\right|$$ is locally
uniformly continuous.\\

\end{lem}
\paragraph{\bf Proof.}
Let $b_0\in]0,1[$ and $m>0$. It is sufficient  to establish,
\begin{equation}\label{continuity}
\left(\sup_{|\gamma|\leq m,\, b\leq b_0} \left|\mathcal
S(\gamma+\tau,b+\delta)-\mathcal S(\gamma,b)\right|\right)
\to0\quad{\rm as~}(\tau,\delta)\to0\,.
\end{equation}
Let $\tau_1=1-b_0>0$. By monotonicity, for all
$\tau,\delta\in[-\tau_1,\tau_1]$, the following holds
$$\{\xi\in\mathbb R~:~\mu_j(\gamma+\tau,\xi)<b+\delta\}\subset
\{\xi\in\mathbb
R~:~\mu_1(-m-\tau_1,\xi)<b+\tau_1\}\,,\quad\forall~j\in\mathbb
N\,.$$ Therefore, we may find a constant $M>0$ depending only on $m$
and $b_0$
 such that
\begin{equation}\label{uniform-bound}
\{\xi\in\mathbb R~:~\mu_j(\gamma+\tau,\xi)<b+\delta\}\subset[-M,M]\,,\quad
\forall~\tau,\delta\in[-\tau_1,\tau_1],~\forall~j\in\mathbb N\,.
\end{equation}
So defining $\xi_j(M)$ as in the proof of Lemma~\ref{justification},
i.e.
$$\forall~\xi\in[-M,M],~\forall~\tau\in[-\tau_1,\tau_1],\quad
\mu_j(\gamma+\tau,\xi)\geq\mu_j(-m-\tau_1,\xi_j(M))\,,$$
we get as in (\ref{claim-j}):
$$\lim_{j\to\infty}\mu_j(-m-\tau_1,\xi_j(M))=\infty\,.$$
Hence, we may find $j_0\geq1$ depending only on $m$ and $b_0$
such that
$$\mu_j(-m-1,\xi_j(M))\geq b_0+2\tau_1\quad \forall~j\geq
j_0\,,$$
and consequently, for $|\tau|\leq\tau_1$, $|\delta|\leq\tau_1$, we get
$$\sum_{j=1}^\infty\left|\{\xi\in\mathbb
R~:~\mu_j(\gamma+\tau,\xi)<b+\delta\}\right|=\sum_{j=1}^{j_0}\left|\{\xi\in\mathbb
R~:~\mu_j(\gamma+\tau,\xi)<b+\delta\}\right|\,.$$ Therefore, we deal
only with a finite sum of $j_0$ terms, $j_0$ being independent from
$\tau$, $\delta$, $\gamma$ and $b$. So given $k\in\{1,\cdots, j_0\}$
and setting $\mathcal S_k(\gamma,b)=\left| \{\xi\in\mathbb
R~:~\mu_k(\gamma,\xi)<b\}\right|$, it is sufficient to show that
\begin{equation}\label{continuity-2}
\lim_{\substack{(\tau,\delta)\to0\\
|\tau|+|\delta|\leq\tau_1}} \left(\sup_{|\gamma|\leq m,\,b\leq
b_0}\left|\mathcal S_k(\gamma+\tau,b+\delta)-\mathcal
S_k(\gamma,b)\right|\right)=0\,.
\end{equation}
The above formula is only a direct consequence of
Lemmas~\ref{St-Li-muk} and \ref{reg-xi-k}.\hfill$\Box$\\

\subsection{The model operator on a half-cylinder}\label{Sub-Sec-MO-2}
We treat now the operator
$P_{h,\Omega_S}^{\alpha,\gamma}=-(h\nabla-iA_0)^2$,
where
$\Omega_S$ is the half-cylinder
$$\Omega_S=]0,S[\times]0,\infty[,$$
$S>0$ and $\gamma\in\mathbb R$ are  constants. The magnetic potential $A_0$ is
taken in the canonical way
\begin{equation}\label{JMP07-A0}
A_0(s,t)=(-t,0),\quad \forall~(s,t)\in[ 0,S]\times[0,\infty[.\end{equation}
Functions in the domain of
$P_{h,\Omega_S}^{\alpha,\gamma}$ satisfy the periodic conditions
$$u(0,\cdot)=u(S,\cdot)\quad {\rm on}~\mathbb R_+,$$
and the de\,Gennes boundary condition at $t=0$,
$$h\,\partial_tu\,\big{|}_{t=0}=h^\alpha\gamma u\,\big{|}_{t=0}.$$

We shall from now on use the following notation. For a self-adjoint
operator $T$ and a real number $\lambda<\inf\sigma_{\rm ess}(T)$, we
denote by $N(\lambda,T)$ the number of eigenvalues of $T$ (counted
with multiplicity)  included in
$]-\infty,\lambda]$.

\begin{lem}\label{JMP07-lem-N-MO}
For each $M>0$,  there exists a constant $C>0$ such that, for all
$b_0\in]-\infty,1[$, $\alpha\geq\frac12$, $S>0$, $\gamma\in[-M,M]$
and $h\in]0,1[$, we have,
\begin{align*}
\begin{split}
&\left|N\left(b_0h,P_{h,\Omega_S}^{\alpha,\gamma}\right)
-\frac{S\,h^{-1/2}}{2\pi}\left(\sum_{j=1}^\infty
\left|\{\xi\in\mathbb R~:~\mu_j(h^{\alpha-1/2}\gamma,\xi)\leq
b_0\}\right|\right)\right|\\
&\leq C.
\end{split}
\end{align*}
\end{lem}
\paragraph{\bf Proof.}
By separation of variables (cf. \cite{ReSi}) and a scaling we may
decompose $P_{h,\Omega_S}^{\alpha,\gamma}$ as a direct sum:
$$\bigoplus_{n\in\mathbb Z}\,h\left(-\frac{\md ^2}{\md \tau^2}+
(2\pi nh^{1/2}S^{-1}+\tau)^2\right)\quad {\rm in}~
\bigoplus_{n\in\mathbb Z} L^2(\mathbb R_+),$$
with the boundary condition  $u'(0)=h^{\alpha-1/2}\gamma\,u(0)$
at $\tau=0$.\\
Consequently we obtain:
\begin{equation}\label{conseq-D1}
\sigma\left(P_{h,\Omega_S}^{\alpha,\gamma}\right)
=\bigcup_{j\in\mathbb N}\left\{h\mu_j(h^{\alpha-1/2}\gamma,2\pi
h^{1/2}S^{-1}n)~:~n\in\mathbb N\right\},\end{equation}
and
each eigenvalue is of multiplicity $1$.\\
Thus, putting
$$f_j(\xi)=\mathbf 1_{\{\xi\in\mathbb
  R~:~\mu_j(h^{\alpha-1/2}\gamma,\xi)<b_0\}}(\xi)\,,\quad
\xi_n=2\pi h^{1/2}S^{-1}n\,,$$
we obtain
\begin{eqnarray*}
N\left(b_0h,P_{h,\Omega_S}^{\alpha,\gamma}\right)&=&
{\rm Card}\left(\sigma
\left(P_{h,\Omega_S}^{\alpha,\gamma}\right)\cap\,]0,b_0h]\right)\\
&=&\sum_{j=1}^\infty
{\rm Card}\left\{n\in\mathbb N~:~\mu_j(h^{\alpha-1/2}\gamma,2\pi
h^{1/2}S^{-1}n)\leq b_0\right\}\\
&=&\sum_{j=1}^\infty \sum_{n\in\mathbb
  Z}f_j(\xi_n)=\sum_{j=1}^{j_0}\sum_{n\in\mathbb Z}f_j(\xi_n)\,.
\end{eqnarray*}
Notice that the last step is due to Lemma~\ref{justification} (and its proof) which
yields the existence of $j_0\in\mathbb N$, depending only on $\gamma$,
$h$ and $b_0$,
such that $f_j\equiv 0$ for $j\geq j_0$.\\
Now, by definition of $\xi_n$,
$$\sum_{n\in\mathbb
  Z}f_j(\xi_n)=\frac{S\,h^{-1/2}}{2\pi}\sum_{n\in\mathbb
  Z}\left(\xi_{n+1}-\xi_n\right)f_j(\xi_n)\,,$$
and we can verify easily  the following estimate (thanks in
  particular to Lemma~\ref{St-Li-muk}),
$$
-2\pi h^{1/2}S^{-1}+\int_{\mathbb R}f_j(\xi)\,\md\xi\leq
\sum_{n\in\mathbb Z}(\xi_{n+1}-\xi_n)f_j(\xi_n)\leq \int_{\mathbb R}f_j(\xi)\,\md \xi\\
+2\pi h^{1/2}S^{-1}\,.
$$
Therefore, we conclude upon noticing that, by  the definition of the function
$f_j$,
$$\int_{\mathbb R}f_j(\xi)\,\md \xi=\left|\{\xi\in\mathbb R~
:~\mu_j(h^{\alpha-1/2}\gamma,\xi)<b_0\}\right|\,.$$
\hfill$\Box$\\

\subsection{The model operator on a
Dirichlet strip}\label{Sub-Sec-MO-3}
We consider now the operator
$P_{h,\Omega_{S,T}}^{\alpha,\gamma}=-(h\nabla-iA_0)^2$, where
$\Omega_{S,T}$ is the strip
$$\Omega_{S,T}=]0,S[\times]0,T[,$$
$S,T>0$ and $\gamma\in\mathbb R$ are constants. The magnetic potential
$A_0$ was defined in (\ref{JMP07-A0}).\\
Functions in the domain of $P_{h,\Omega_{S,T}}^{\alpha,\gamma}$
satisfy the de\,Gennes
condition
$h\partial_tu=h^\alpha\gamma\,u$ at $t=0$ and Dirichlet
condition on the other sides of the boundary.\\

The next lemma  gives a comparison between the counting function of
$P_{h,\Omega_{S,T}}^{\alpha,\gamma}$ and that of
$P_{h,\Omega_S}^{\alpha,\gamma}$.

\begin{lem}\label{JMP07-lem-N-MO'}
There exists a constant $c>0$ such that,
$$\forall~S>0,~\forall~T>0,\quad \forall~\gamma\in\mathbb R,\quad\forall
~\delta\in]0,S/2],\quad\forall~b_0\in]-\infty,1[,$$
we have,
\begin{eqnarray*}
\frac12N\left(b_0h-c\,h^2(\delta^{-2}+T^{-2}),\,
P_{h,\Omega_{2(S-\delta)}}^{\alpha,\gamma}\right)&\leq&
N\left(b_0h,P_{h,\Omega_{S,T}}^{\alpha,\gamma}\right)\\
&\leq&
N\left(b_0h,P_{h,\Omega_S}^{\alpha,\gamma}\right).
\end{eqnarray*}
\end{lem}
\paragraph{\bf Proof.}
Since the extension by zero of a  function in the  form domain of $P_{h,\Omega_{S,T}}^{\alpha,\gamma}$ is
included in that of $P_{h,\Omega_S}^{\alpha,\gamma}$, and the values
of the quadratic forms coincide for such a function, we get the upper
bound of the lemma by a simple application of the variational
principle.\\
We turn now to the lower bound. The argument  is like the one used in
\cite{CV, Fr} but we   explain it because it
illustrates in a simple case the arguments of this paper.\\
Let us introduce two  partitions of unity $(\varphi_i^\delta)$ and
$(\psi^T_j)$ such that:
$$\left(\varphi^\delta_1\right)^2+\left(\varphi_2^\delta\right)^2=1~
{\rm  in}~[0,2(S-\delta)],\quad
\left(\psi_0^T\right)^2+\left(\psi_1^T\right)^2=1~{\rm in}~\mathbb
R_+,$$
$$\left\{\begin{array}{c}
{\rm supp}\,\varphi_1^\delta\subset[0,S]\\
{\rm supp}\,\varphi_2^\delta\subset[S-\delta,2(S-\delta)]\\
\displaystyle\sum_{i=1}^2\left|\left(\varphi_i^\delta\right)'\right|^2\leq
c\,\delta^{-2}
\end{array}\right\}~\&~
\left\{\begin{array}{c}
{\rm supp}\,\psi_0^T\subset[T/2,\infty[\\
{\rm supp}\,\psi_1^T\subset[0,T]\\
\displaystyle\sum_{i=0}^1\left|\left(\psi_i^T\right)'\right|^2\leq
c\,T^{-2}
\end{array}\right\},$$
where $c>0$ is a constant independent from $S,T$ and $\delta$.\\
Upon putting
$$\chi_i^{\delta,T}(s,t)=\varphi_i^\delta(s)\,\psi_1^T(t)\quad
(i=1,2),\quad
\chi_0^{\delta,T}(s,t)=\psi^T_0(t),$$
we get a partition of unity of $\Omega_{2(S-\delta)}=
]0,2(S-\delta)[\times\mathbb R_+$.\\
Let us take a function $u$ in the form domain of
$P_{h,\Omega_{2(S-\delta)}}^{\alpha,\gamma}$. Then, by the IMS
decomposition formula:
\begin{eqnarray*}
&&\hskip-1cm\int_{\Omega_{2(S-\delta)}}|(h\nabla-iA_0)u|^2\md x\\
&&=\sum_{i=0}^2\int_{\Omega_{2(S-\delta)}}
|(h\nabla-iA_0)\,\chi_i^{\delta,T}u|^2\md x-h^2\sum_{i=0}^2\left\|\,
|\nabla\chi_i^{\delta,T}|u\right\|^2_{L^2(\Omega_{S,T})}\\
&&\geq \sum_{i=0}^2\int_{\Omega_{2(S-\delta)}}
|(h\nabla-iA_0)\,\chi_i^{\delta,T}u|^2\md
x-c(\delta^{-2}+T^{-2})h^2\|u\|^2_{L^2(\Omega_{S,T})}.
\end{eqnarray*}
Notice that $\chi_1^{\delta,T}u$ is in the form domain of
$P_{h,\Omega_{S,T}}^{\alpha,\gamma}$ and $\chi_2^{\delta,T}u$ is in
that of $P_{h,]S-\delta,2(S-\delta)[\times]0,T[}^{\alpha,\gamma}$
(this last operator, thanks to translational invariance with respect
to $s$, is unitary equivalent to
$P_{h,\Omega_{S,T}}^{\alpha,\gamma}$). Also, $\chi_0^{\delta,T}u$ is
in the form domain of the $P_{h,\Omega_{2(S-\delta)}}^D$, the
Dirichlet realization of $-(h\nabla-iA_0)^2$ in
$\Omega_{2(S-\delta)}$.\\
Since the form domain of
$P_{h,\Omega_{2(S-\delta)}}^{\alpha,\gamma}$ can be viewed
in\footnote{For an operator $A$, $FD(A)$ denotes its form domain.}
$$FD\left(P_{h,\Omega_{S,T}}^{\alpha,\gamma}\right)\oplus FD\left(P_{h,]S-\delta,2(S-\delta)[\times]0,T[}^{\alpha,\gamma}\right)\oplus
FD\left(P_{h,\Omega_{2(S-\delta)}}^D\right)$$ via the isometry
$u\mapsto
\left(\chi_1^{\delta,T}u,\chi_2^{\delta,T}u,\chi_0^{\delta,T}u\right)$,
we get upon applying the variational principle (see
\cite[Section~XII.15]{ReSi}),
\begin{eqnarray}\label{JMP07-lem-MO'-Proof}
&&\hskip-1cm
N\left(b_0h-c\,h^2(\delta^{-2}+T^{-2}),\,P_{h,\Omega_{2(S-\delta)}}^{\alpha,\gamma}\right)\\
&&\hskip1cm\leq 2 N\left(b_0h,\,P_{h,\Omega_{2(S-\delta)}}^
{\alpha,\gamma}\right)+N\left(b_0h,\,P_{h,\Omega_{2(S-\delta)}}^D\right).\nonumber\end{eqnarray}
Notice that, by (\ref{JMP07-Simon}), the operator
$P_{h,\Omega_{2(S-\delta)}}^D$ has no spectrum below $b_0h$ when
$b_0<1$. Hence,
$$N\left(b_0h,\,P_{h,\Omega_{2(S-\delta)}}^D\right)=0.$$
Coming back to (\ref{JMP07-lem-MO'-Proof}),  we get the lower bound
stated in the lemma.\hfill$\Box$\\

\section{Proof of Theorem~\ref{JMP07thm1}}\label{proofThm1}
We come back to the case of a general smooth domain $\Omega$ whose
boundary is compact. We introduce the following quadratic forms:
\begin{eqnarray}\label{VI-qform}
&&q_{h,\Omega}^{\alpha,\gamma}(u)=\int_{\Omega}|(h\nabla-iA)u|^2\,\md x+
h^{1+\alpha}\int_{\partial\Omega}\gamma(s)\,|u(s)|^2\,\md s,\\
&&q_{h,\Omega}(u)=\int_{\Omega}|(h\nabla-iA)u|^2\,\md
x,\label{VI-qform0}
\end{eqnarray}
defined for functions in the magnetic Sobolev space:
\begin{equation}\label{VI-H-A1}
H_{h,A}^1(\Omega)=\{u\in L^2(\Omega)~:~(h\nabla -iA)u\in
L^2(\Omega)\}.\end{equation}
Here, we recall that $A\in C^2(\overline{\Omega})$ is such that
$${\rm curl}\,A=1\quad {\rm in}~\overline{\Omega}.$$
We shall recall in the appendix a standard coordinate transformation valid in a
sufficiently thin neighborhood of the boundary:
$$\Phi_{t_0}:~\overline{\Omega(t_0)}\ni x\mapsto (s(x),t(x))\in
[0,
|\partial\Omega|[\,\times\,[0,t_0],$$
where for $t_0>0$, $\Omega(t_0)$ is the tubular neighborhood of
$\partial\Omega$:
$$\Omega(t_0)=\{x\in\Omega~:~{\rm dist}(x,\partial\Omega)<t_0\}.$$
Let us mention that $t(x)={\rm dist}(x,\partial\Omega)$
measures the distance to the boundary and
$s(x)$ measures the curvilinear distance in $\partial\Omega$.\\
Using the coordinate transformation $\Phi_{t_0}$, we associate to any
function $u\in L^2(\Omega)$, a function $\widetilde u$ defined in $
[0,
|\partial\Omega|[\,\times\,[0,t_0]$ by,
\begin{equation}\label{VI-utilde}
\widetilde u(s,t)=u(\Phi_{t_0}^{-1}(s,t)).
\end{equation}

The next lemma
states  a standard approximation of the quadratic form
$q_{h,\Omega}^{\alpha,\gamma}(u)$ by the canonical one
in the half-plane, provided
that the function $u$ is supported near the boundary.

\begin{lem}\label{VI-FrLem3.5}
There exists a constant $C>0$, and for all
$$S_1\in[0,|\partial\Omega|[,\quad
S_2\in]S_1,|\partial\Omega|[,$$
there exists a function $\phi\in C^2([S_1,S_2]\times[0,t_0];\mathbb
R)$ such that, for all
$$\widetilde S_1\in[S_1,S_2],\quad T\in
  ]0,t_0[,\quad
\varepsilon\in[CT,Ct_0],$$
and for all $u\in H_{h,A}^1(\Omega)$ satisfying
$${\rm supp}\,\widetilde u\subset [S_1,S_2]\times[0,T],$$
one has the following estimate,
\begin{eqnarray*}
&&\hskip-0.5cm(1-\varepsilon)
q_{h,\Omega_1}^{\alpha,\widehat\gamma_1}\left(e^{i\phi/h}\widetilde u\right)
-C\varepsilon^{-1}\left(((S_2-S_1)^2+T^2)^2+h^2\right)
\|\widetilde u\|_{L^2(\Omega_1)}^2\\
&&\leq\, q_{h,\Omega}^{\alpha,\gamma}(u)
\leq\,(1+\varepsilon)
q_{h,\Omega_1}^{\alpha,\widetilde\gamma_1}\left(e^{i\phi/h}\widetilde u\right)
+C\varepsilon^{-1}\left(((S_2-S_1)^2+T^2)^2+h^2\right)
\|\widetilde u\|_{L^2(\Omega_1)}^2.
\end{eqnarray*}
Here $\Omega_1=[S_1,S_2]\times[0,T]$, $\widetilde\gamma_1=
\displaystyle\frac{\gamma(\widetilde S_1)+C(S_2-S_1)}{1+\varepsilon}$,
$\widehat\gamma_1=
\displaystyle\frac{\gamma(\widetilde S_1)-C(S_2-S_1)}{1-\varepsilon}$,
and the function $\widetilde u$ is associated to $u$ by (\ref{VI-utilde}).
\end{lem}
\paragraph{\bf Proof.}
Notice that when $\gamma\equiv0$, the result  follows from
\cite[Lemma~3.5]{Fr}, which reads explicitly in the form:
\begin{eqnarray*}
&&\hskip-1cm\left|q_{h,\Omega}^{\alpha,0}(u)-
\int_{[S_1,S_2]\times[0,T]}|(h\nabla-iA_0) e^{i\phi/h}\widetilde
u|^2\md s\md t\right|\\
&&\leq \varepsilon
\int_{[S_1,S_2]\times[0,T]}|(h\nabla-iA_0)e^{i\phi/h}\widetilde
u|^2\md s\md
t\\
&&\hskip0.5cm+C\varepsilon^{-1}\left((T^2+(S_2-S_1)^2)^2+h^2\right)
\|\widetilde u\|_{L^2}^2.
\end{eqnarray*}
Since $\widetilde u$ restricted to the boundary is supported in
$[S_1,S_2]$,
we get as an immediate consequence the following two-sided estimate
for non-zero $\gamma$:
\begin{align}\label{eq:Sref1}
\begin{split}&\hskip-0.1cm(1+\varepsilon)^{-1}q_{h,\Omega}^{\alpha,\gamma}(u)\\
& \leq\int_{[S_1,S_2]\times[0,T]}|
(h\nabla-iA_0)e^{i\phi/h}\widetilde u|^2\md s\md
t+\frac{h^{1+\alpha}}{1+\varepsilon}\int_{[S_1,S_2]}\gamma(s)
|\widetilde u(s,0)|^2\md s\\
&\hskip0.5cm+C(1+\varepsilon)^{-1}\varepsilon^{-1}\left((T^2+(S_2-S_1)^2)^2+h^2\right)
\|\widetilde u\|_{L^2}^2\,,
\end{split}
\end{align}
and
\begin{align}\label{eq:Sref2}
\begin{split}
&\hskip-0.1cm(1-\varepsilon)^{-1}
q_{h,\Omega}^{\alpha,\gamma}(u)\\
&\geq\int_{[S_1,S_2]\times[0,T]}|
(h\nabla-iA_0)e^{i\phi/h}\widetilde u|^2\md s\md
t+\frac{h^{1+\alpha}}{1-\varepsilon}\int_{[S_1,S_2]}\gamma(s)
|\widetilde u(s,0)|^2\md s\\
&\hskip0.5cm-C(1-\varepsilon)^{-1}\varepsilon^{-1}\left((T^2+(S_2-S_1)^2)^2+h^2\right)
\|\widetilde u\|_{L^2}^2\,.
\end{split}
\end{align}
The idea now is to approximate $\gamma$ by a constant in a simple
manner without needing an estimate of the boundary integral.
Actually, Taylor's formula applied to the function $\gamma$ near
$\widetilde S_1$ leads to the estimate
$$|
\gamma(s)-\gamma(\widetilde
S_1)|\leq C(S_2-S_1)\,,\quad\forall~s\in[S_1,S_2]\,,$$
where the constant $C>0$ is possibly replaced by a larger one.\\
Having this estimate in hand, we get:
\begin{eqnarray*}
&&\hskip-0.5cm
\left|\int_{[S_1,S_2]}\gamma(s)
|\widetilde u(s,0)|^2\md s-
\gamma(\widetilde S_1)\int_{[S_1,S_2]}|\widetilde
u(s,0)|^2\md s\right|\\
&&\leq C(S_2-S_1)\int_{[S_1,S_2]}|\widetilde u(s,0)|^2\,\md s\,.
\end{eqnarray*}
Recalling the definition of $\widetilde\gamma$ and $\widehat\gamma$
in Lemma~\ref{VI-FrLem3.5} (they actually depend on $\varepsilon$,
$S_1$, $S_2$ and hence account to all possible errors), we infer
directly from the previous estimate, \begin{align}\label{eq:Sref3}
\begin{split}
(1-\varepsilon)\widehat \gamma_1 \int_{[S_1,S_2]} |\widetilde
u(s,0)|^2\md s&\leq
\int_{[S_1,S_2]}\gamma(s) |\widetilde u(s,0)|^2\md s\\
&\leq (1+\varepsilon)\widetilde\gamma_1\int_{[S_1,S_2]} |\widetilde
u(s,0)|^2\md s\,.\end{split}\end{align} Substituting the lower and
upper bound of (\ref{eq:Sref3}) in (\ref{eq:Sref1}) and
(\ref{eq:Sref2}) respectively, and recalling the hypothesis that
$\widetilde u(s,0)$
is supported in $[S_1,S_2]$, we obtain the desired estimates of the lemma.\hfill$\Box$\\

We shall divide a thin neighborhood of $\partial\Omega$ into
many small sub-domains, and in each sub-domain, we shall apply
Lemma~\ref{VI-FrLem3.5} to approximate the quadratic form. This will
yield a two-sided estimate of
$N(\lambda, P_{h,\Omega}^{\alpha,\gamma})$ in terms of
the spectral counting functions of model operators on
half-cylinders.\\
Let us put
$$N=[h^{-{3/8}}],$$
the greatest  positive integer below $h^{-{3/8}}$.\\
Let
$$S=\frac{|\partial\Omega|}N,\quad s_n=nS,~~n\in\{0,1,\dots,N\},$$
and we emphasize that these quantities depend on $h$. We put further
$$\Omega_S=]0,S[\times\mathbb R.$$

With these notations, the proof of Theorem~\ref{JMP07thm1} is given by
the following proposition.

\begin{prop}\label{VI-FrProp3.6}
Let $b_0\in]-\infty,1[$. There exist constants $C,h_0>0$ such that
for all
$$h\in]0,h_0],\quad \delta\in]0,S/2],\quad \widetilde
S_n\in[s_{n-1},s_n],\quad
n\in\{1,2,\dots,N\},$$
one has the  following estimate on the spectral counting function,
\begin{eqnarray*}
\frac12\sum_{n=1}^N
N\left(hb_0-Ch^2\delta^{-2},P_{h,\Omega_{2(S-\delta)}}^{\alpha,\widetilde\gamma_n}\right)
&\leq&
N\left(hb_0,P_{h,\Omega}^{\alpha,\gamma}\right)\\
&\leq&
\sum_{n=1}^N
N\left(hb_0+Ch^2\delta^{-2},
P_{h,\Omega_{S+2\delta}}^{\alpha,\widehat\gamma_n}\right).
\end{eqnarray*}
Here $\widetilde\gamma_n=\displaystyle\frac{\gamma(\widetilde
  S_n)+CS}{1+h^{1/4}}$
and $\widehat\gamma_n=\displaystyle\frac{\gamma(\widetilde S_n)-CS}{1-h^{1/4}}$.
\end{prop}

Before proving Proposition~\ref{VI-FrProp3.6}, let us see how it serves
for obtaining  the conclusion of Theorem~\ref{JMP07thm1}.

\paragraph{\bf Proof of Theorem~\ref{JMP07thm1}}
We keep the notation introduced for the statement of
Proposition~\ref{VI-FrProp3.6}. The proof is in two steps.\\
{\it Step~1.}\\
Let us establish the asymptotic formula (as $h\to0$):
\begin{eqnarray}\label{JMP07-proofthm1}
&&h^{1/2}N\left(b_0h,P_{h,\Omega}^{\alpha,\gamma}\right)
=\\
&&\hskip0.7cm\frac1{2\pi}\int_{\partial\Omega}\sum_{j=1}^\infty\big{|}
\{\xi\in\mathbb
R~:~\mu_j(h^{\alpha-1/2}\gamma(s),\xi)<b_0\}\big{|}\md
s+o(1),\nonumber
\end{eqnarray}
where $b_0\in]-\infty,1[$.\\
From the lower bound in Proposition~\ref{VI-FrProp3.6}, we get upon applying
  Lemma~\ref{JMP07-lem-N-MO}, a constant $\widetilde C>0$ such that
\begin{eqnarray*}
&&\hskip-0.5cmh^{1/2}N\left(b_0h,P_{h,\Omega}^{\alpha,\gamma}\right)
\geq\\
&&\hskip0.5cm\frac1{2\pi}\sum_{n=1}^N(S-\delta)\sum_{j=1}^\infty
\left|\{\xi\in\mathbb R~:~\mu_j(h^{\alpha-1/2}\widetilde
\gamma_n,\xi)<b_0-Ch\delta^{-2}\}\right|-\widetilde C h^{1/2}.
\end{eqnarray*}
Since $\alpha\geq1/2$ and $\partial\Omega$ is bounded, it is a
result of Lemmas~\ref{justification} and \ref{lem-referee} that
there exist constants $C>0$ and $h_0>0$ together with a function
$]0,h_0]\ni h\mapsto\epsilon(h)$ tending to $0$ as $h\to 0$ such that,
for all $h\in]0,h_0]$ and
$n\in\{1,2,\dots,N\}$, we have (provided that $h\delta^{-2}$ is
sufficiently small),
\begin{eqnarray*}
&&
\mathcal S\left(h^{\alpha-1/2}\widetilde\gamma_n,b_0\right)\leq C\,,\\
&&\left|\mathcal
S\left(h^{\alpha-1/2}\widetilde\gamma_n,b_0-Ch\delta^{-2}\right)
-\mathcal S\left(h^{\alpha-1/2}\gamma_n,b_0\right)\right|\leq\epsilon(h)\,,
\end{eqnarray*}
where the function $\mathcal S$ is introduced in \ref{lem-referee}, and
$\gamma_n=(1+h^{1/4})\widetilde\gamma_n-CS=\gamma(\widetilde
S_n)$\,.\\
Recalling that  $S=|\partial\Omega|/N$ and $N=[h^{-3/8}]$, we get
upon choosing $\delta=1/N^2$,
\begin{eqnarray*}
&&\hskip-0.5cm h^{1/2}N\left(b_0h,P_{h,\Omega}^{\alpha,\gamma}\right)\\
&&\geq \frac1{2\pi}\left(\sum_{n=1}^N\sum_{j=1}^\infty
S\left|\{\xi\in\mathbb R~:~\mu_j(h^{\alpha-1/2}
\gamma_n,\xi)<b_0\}\right|\right)-\epsilon(h)-Ch^{3/8},\end{eqnarray*}
where the
leading order term on the right hand side is a Riemann sum. We get
then the following lower bound,
\begin{eqnarray*}&&\hskip-0.5cm
h^{1/2}N\left(b_0h,P_{h,\Omega}^{\alpha,\gamma}\right)
\geq\\
&&\hskip1cm\frac1{2\pi}\int_{\partial\Omega}\sum_{j=1}^\infty
\left|\{\xi\in\mathbb R~:~\mu_j(h^{\alpha-1/2}
\gamma(s),\xi)<b_0\}\right|\,\md s+o(1).\end{eqnarray*} This is the
lower bound in (\ref{JMP07-proofthm1}). In a similar manner
we obtain an upper bound.\\
{\it Step~2.}\\
If $\alpha=1/2$, the asymptotic formula (\ref{JMP07-proofthm1}) is
just the conclusion of Theorem~\ref{JMP07thm1}.\\
We turn to the case when $\alpha>1/2$. Again, it results from
Lemma~\ref{lem-referee} the existence of a constant $h_0$ and a
function $]0,h_0]\ni h\mapsto\epsilon_1(h)$ tending to $0$ as
$h\to0$ such that for all $h\in]0,h_0]$ and $s\in\partial\Omega$,
$$\left|\mathcal
S\left(h^{\alpha-1/2}\gamma(s),b_0\right)
-\left(\sum_{j=1}^\infty|\{\xi\in\mathbb R~:~\mu_j(\xi)<b_0|\right)\right|
\leq \epsilon_1(h),$$
where $\mu_j(\xi)=\mu_j(0,\xi)$. Moreover, by
Lemma~\ref{RuFr-Lemm-kach}, it holds that
$$\sum_{j=2}^\infty\left|\{\xi\in\mathbb
R~:~\mu_j(\xi)<b_0\}\right|=0\,.$$ Now  we
can infer from  (\ref{JMP07-proofthm1}) the
asymptotic formula announced in (\ref{JMP07thm11}).\hfill$\Box$\\

\paragraph{\bf Proof of Proposition~\ref{VI-FrProp3.6}.}
Let us establish the lower bound.
Let $P_{N,h,\Omega}^{\alpha,\gamma}$
be the restriction of the operator $P_{h,\Omega}^{\alpha,\gamma}$ for
functions $u$ in $D(P_{h,\Omega}^{\alpha,\gamma})$ that vanishes on the
set
$$\{x\in\Omega~:~t(x)\geq T\}\cup
\bigcup_{n=1}^N\{x\in\Omega~:~0\leq t(x)\leq T, s(x)=s_n\},$$
where $T>$ is to be specified later. The important remark is that the
spectrum of $P_{h,\Omega}^{\alpha,\gamma}$ is below that of
$P_{N,h,\Omega}^{\alpha,\gamma}$.\\
Let us take a function $u$ in the form of domain of
$P_{N,h,\Omega}^{\alpha,\gamma}$. Applying Lemma~\ref{VI-FrLem3.5}
with $T=h^{3/8}$ and $\varepsilon=h^{1/4}$, we get the estimate,
\begin{eqnarray*}
q_{h,\Omega}^{\alpha,\gamma}(u)\leq
(1+h^{1/4})\sum_{n=1}^N\left(
q_{h,\Omega_n}^{\alpha,\widetilde\gamma_n}\left(e^{-i\phi_n/h}\,\widetilde
u\right)+Ch^{5/4}\left\|e^{-i\phi_n/h}\,\widetilde
u\right\|_{L^2(\Omega_n)}^2\right),
\end{eqnarray*}
where $\Omega_n=]s_{n-1},s_n[\,\times\,]0,T[$ and
$\widetilde\gamma_n= \frac{\gamma(\widetilde S_n)+CS}{1+h^{1/4}}$.
Then, by the variational principle, we obtain (recall that the
    spectrum of $P_{h,\Omega}^{\alpha,\gamma}$ is below that of
    $P_{N,h,\Omega}^{\alpha,\gamma}$),
$$N\left(\lambda,\,P_{h,\Omega}^{\alpha,\gamma}\right)\geq
\sum_{n=1}^NN\left(\frac{\lambda-Ch^{5/4}}{1+h^{1/4}},
P_{h,\Omega_{S,T}}^{\alpha,\widetilde \gamma_n}\right),$$
where $\Omega_{S,T}=]0,S[\times]0,T[$. Applying
  Lemma~\ref{JMP07-lem-N-MO'}, this is sufficient to conclude the lower
  bound announced in Proposition~\ref{VI-FrProp3.6}.\\
The upper bound is derived by introducing a partition of unity
attached to the sub-domains $\Omega_n$ and by using the IMS
decomposition  formula.
The analysis is similar to that
presented for the lower bound above and also to that in the proof of
Lemma~\ref{JMP07-lem-N-MO'}, so we omit the proof. For the details,
we refer to \cite[Proposition~3.6]{Fr}.\hfill$\Box$\\

\section{Curvature effects}\label{CurvEffects}
\subsection{A family of ordinary differential operators on a weighted $L^2$
space}\label{CurvEffects-MO}\ \\
A finer approximation of  the quadratic form (\ref{VI-qform})
leads to the analysis of a family of ordinary differential
operators on a weighted $L^2$ space that takes into account the
curvature effects of the boundary. We shall recall in this section
the main results for the lowest eigenvalue problem concerning this
family of operators. These results were obtained in \cite{HeMo}
for the Neumann problem and then generalized in
\cite{Kach1} for situations involving de\,Gennes' boundary condition.\\
Let us introduce, for technical reasons that will be clarified
later, a positive parameter $\delta\in]\frac14,\frac12[$. Let us
also consider parameters $h>0$ and $\beta\in\mathbb R$ such that
$$|\beta| h^\delta<\frac13.$$
We define the family of quadratic forms (indexed by $\xi\in\mathbb
R$)
\begin{eqnarray}\label{Th-qf}
q_{h,\beta,\xi}^{\alpha,\eta}(u)&=&\int_0^{h^{\delta-1/2}}\left[
|u'(t)|^2+(1+2\beta h^{1/2}t)\left|\left(t-\xi-\beta
h^{1/2}\frac{t^2}2\right)u(t)\right|^2\right]\nonumber\\
&&\hskip4cm\times(1-\beta
h^{1/2}t)\,dt+h^{\alpha-1/2}\eta|u(0)|^2,
\end{eqnarray}
defined for functions $u$ in the space~:
\begin{equation}\label{Th-df}
D\left(q_{h,\beta,\xi}^{\alpha,\eta}\right)=\left\{u\in
H^1\left(]0,h^{\delta-1/2}[\right)~:~u\left(h^{\delta-1/2}\right)=0\right\}.
\end{equation}
Let us denote by $H_{h,\beta,\xi}^{\alpha,\eta}$ the self-adjoint
realization associated to the quadratic form (\ref{Th-qf}) by
Friedrich's' theorem. Let us denote also by
$\left(\mu_j\left(H_{h,\beta,\xi}^{\alpha,\eta}\right)\right)_j$
the increasing sequence of eigenvalues of
$H_{h,\beta,\xi}^{\alpha,\eta}$.\\
For each $\alpha\geq\frac12$ and $\eta\in\mathbb R$ we introduce the
positive numbers~:
\begin{align}\label{CE-d2-d3}
\begin{split}
&d_2\left(\frac12,\eta\right)=\xi(\eta)\Theta'(\eta),\quad
d_2(\alpha,\eta)=\xi_0\Theta'(0)\quad\left(\alpha>\frac12\right),\\
&d_3\left(\frac12,\eta\right)=
\frac13(\eta\xi(\eta)+1)^2\Theta'(\eta),\quad
d_3(\alpha,\eta)=\frac13\Theta'(0)\quad\left(\alpha>\frac12\right).
\end{split}\end{align}

The result concerning
the lowest eigenvalue of
$H_{h,\beta,\xi}^{\alpha,\eta}$ in the next theorem has been proved in
\cite{Kach1}.

\begin{thm}\label{Kachmar-LmV.8-V.9}
Suppose that $\delta\in]\frac14,\frac12[$ and $\alpha\geq\frac12$.
Let
$$\widetilde \eta=h^{\alpha-1/2}\eta\,,\quad
\rho_0=\left\{\begin{array}{l}
\delta-\frac14,~{\rm if~}\alpha=\frac12\,\\
\min(\delta-\frac14,\alpha-\frac12),~{\rm if~}\alpha>\frac12
\end{array}\right.\,,\quad{\rm and~}
0<\rho<\rho_0\,.
$$
For each $M>0$  and $\zeta_0>0$,
there exists a constant $\zeta_1>0$, and for each $\zeta\geq\zeta_1$,
there exist positive constants $C,h_0$ and a function
$]0,h_0]\ni h\mapsto\epsilon(h)\in\mathbb R_+$
with $\lim_{h\to0}\epsilon(h)=0$ such that,
$$\forall~\eta,\beta\in]-M,M[,\quad
\forall~h\in]0,h_0],$$
the following assertions hold~:
\begin{itemize}
\item If $|\xi-\xi\left(\widetilde\eta\right)|\leq\zeta h^{\rho}$,
then
\begin{eqnarray}\label{EqappmuG}
&&\left|\mu_1\left(H_{h,\beta,\xi}^{\alpha,\eta,D}\right)-\left\{
\Theta\left(\widetilde\eta\right)
+d_2(\alpha,\eta)(\xi-\xi(\widetilde\eta))^2-
d_3(\alpha,\eta)\beta h^{1/2}\right\}
\right|\\
&&\leq C\left[h^{1/2}|\xi-\xi(\widetilde\eta)|+
h^{\delta+1/2}+h^{1/2}\epsilon(h)\right],\nonumber
\end{eqnarray}
and
\begin{equation}\label{Th-mu2}
\mu_2\left(H_{h,\beta,\xi}^{\alpha,\eta,D}\right)\geq
\Theta\left(\widetilde\eta\right)+\zeta_0h^{2\rho}.
\end{equation}
\item If $\left|\xi-\xi\left(\widetilde\eta\right)\right|\geq
\zeta h^{\rho}$, then
\begin{equation}\label{eq.5.48}
\mu_1\left(H^{\alpha,\eta,D}_{h,\beta,\xi}\right)\geq
\Theta\left(\widetilde\eta\right)+\zeta_0 h^{2\rho}.
\end{equation}
\end{itemize}
Here, the parameters $d_2(\alpha,\eta)$ and $d_3(\alpha,\eta)$ has
been introduced in (\ref{CE-d2-d3}).
\end{thm}
\paragraph{\bf Proof.}
The existence of $\zeta_1$ so that the lower bound (\ref{eq.5.48}) holds for
$|\xi-\xi(\widetilde\eta)|\geq\zeta_1 h^\rho$ has
been established in \cite[Lemma~V.8]{Kach1}. Now, for
$\zeta\geq\zeta_1$, (\ref{eq.5.48}) obviously holds. Under the
hypothesis
$|\xi-\xi(\widetilde\eta)|\leq\zeta h^\rho$, the existence of the
constants
$C$, $h_0$ and the estimate
(\ref{EqappmuG})  have been established in
\cite[Lemma~V.8 \& V.9]{Kach1}. So we only need to establish
(\ref{Th-mu2}).\\
We start with the case $\alpha=\frac12$.
It results from the min-max principle (see \cite{Kach1} or
\cite[Lemma~4.2.1]{KachTh} for details),
\begin{align*}
\left|\mu_2\left(H^{\alpha,\eta,D}_{h,\beta,\xi}\right)-
\mu_2\left(H^{\alpha,\eta,D}_{h,0,\xi}\right)\right|\leq
\widetilde Ch^{2\delta-\frac12}\left(1+
\mu_2\left(H^{\alpha,\eta,D}_{h,0,\xi}\right)\right)\,,
\end{align*}
where the constant $\widetilde C$ depends only on $M$.\\
It results again from the min-max principle,
$$\mu_2\left(H^{\alpha,\eta,D}_{h,0,\xi}\right)\geq
\mu_2\left(\mathcal L[\eta,\xi]\right),$$
where $\mathcal L[\eta,\xi]$ is the operator introduced in (\ref{L-dop}).\\
We get then the following lower bound,
\begin{eqnarray*}
\mu_2\left(H^{\alpha,\eta,D}_{h,\beta,\xi}\right)&\geq&
\left(1-2\widetilde Ch^{2\delta-\frac12}\right)\mu_2\left(\mathcal
L[\eta,\xi]\right)-\widetilde Ch^{2\delta-\frac12}\\
&\geq&\left(1-2\widetilde Ch^{2\delta-\frac12}\right)\Theta_2(\eta)-\widetilde Ch^{2\delta-\frac12}\\
&\geq&\Theta(\eta)+\zeta_0h^{2\rho}.
\end{eqnarray*}
Here, we recall the definition of  $\Theta_2(\eta)$ in
(\ref{Theta2}). Let us also point out that the final conclusion above
follows by Lemma~\ref{HeRu} upon taking $h\in]0,h_0]$ with $h_0$
    chosen so small that $\zeta_0h_0^{2\rho}+2\widetilde
    Ch_0^{2\delta-\frac12}(\Theta_2(\eta)+1)
    <\frac12\left(
\Theta_2(\eta)-\Theta(\eta)\right)$.\\
When $\alpha>\frac12$, the result follows from the above argument upon
using the continuity of our spectral functions with respect to small
perturbations.
\hfill$\Box$\\

\subsection{Spectral function for the model operator on a
  half-cylinder}
\label{CE-MO}\ \\
Let us consider parameters
$$h>0,\quad \delta>0,\quad S>0~{\rm  and}~\beta\in\mathbb R~{\rm s.t.}~
|\beta|h^\delta<\frac12.$$
We denote by $\widetilde L_{h,\beta,S}^{\alpha,\eta}$ the self adjoint operator in
$$L^2\left(]0,S[\times]0,h^\delta[\,;(1-\beta t)\md s\md t\right)$$
associated with the quadratic form
\begin{eqnarray*}
\widetilde Q_{h,\beta,S}^{\alpha,\eta}(u)&=&\int_0^S\int_0^{h^\delta}
\left(|h\partial_tu|^2+(1+2\beta t)\left|\left(h\partial_s+t
-\beta\frac{t^2}2\right)u\right|^2\right)\\
&&\hskip2cm\times(1-\beta t)\,\md s\md
t+h^{1+\alpha}\eta\int_0^S|u(s,0)|^2\,\md s,
\end{eqnarray*}
defined for functions $u$ in the form domain
$$D\left(\widetilde Q_{h,\beta,S}^{\alpha,\eta}\right)=\{u\in H^1\left(
]0,S[\times]0,h^\delta[\right)~:~u(\cdot,h^\delta)=0,~u(0,\cdot)
=u(S,\cdot)\}.$$ We recall again the notation that for a self
adjoint operator $A$ and a number $\lambda<\inf \sigma_{\rm
ess}(A)$, we denote by $N(\lambda ,A)$ the number of eigenvalues of
$A$ (counted with multiplicity) below $\lambda$.

\begin{prop}\label{Fr-prop4.1}
With the notation and hypotheses of Theorem~\ref{Kachmar-LmV.8-V.9}, let
$\zeta_0>0$,
$\widetilde h_0>0$
and $\lambda=\lambda(h)$ such that
\begin{equation}\label{Hyp-lambda}
|\lambda-\Theta(\widetilde\eta)|<\zeta_0 h^{2\rho}, \quad
\forall~h\in]0,\widetilde h_0].\end{equation} Then, for each $M>0$,
there exist constants $C>0$ and $h_0>0$ and a function $]0,h_0]\ni
h\mapsto\epsilon_0(h)\in\mathbb R_+$ with
$\lim_{h\to0}\epsilon_0(h)=0$ such that, for all $h\in]0,h_0]$ and
$S,\eta,\beta\in]-M,M[$,
\begin{align}\label{Conc-N}
\begin{split}
&\left|N\left(h\lambda,\,\widetilde L_{h,\beta,S}^{\alpha,\eta}\right)
-\frac{h^{-1/4}S}{\pi\sqrt{d_2(\alpha,\eta)}}\sqrt{\left(
d_3(\alpha,\eta)\beta+h^{-1/2}[\lambda-\Theta(\widetilde\eta)]\right)_+}\right|
\\
&\leq \frac{h^{-1/4}S}{\pi\sqrt{d_2(\alpha,\eta)}}
\sqrt{\left(
d_3(\alpha,\eta)\beta+h^{-1/2}[\lambda-\Theta(\widetilde\eta)]\right)_+}
\epsilon_0(h).
\end{split}
\end{align}
\end{prop}
\paragraph{\bf Proof.}
By separation of variables and by performing the scaling
$\tau=h^{-1/2}t$, we decompose $\widetilde L_{h,\beta,S}^{\alpha,\eta}$
as a direct sum,
$$\bigoplus_{n\in\mathbb Z}h
\,H_{h,\beta,2\pi n h^{1/2}S^{-1}}^{\alpha,\eta}\quad{\rm in}\quad
\bigoplus_{n\in\mathbb N}L^2\left
(]0,h^{\delta-1/2}[;(1-\beta h^{1/2}t)dt\right).$$
Consequently, by Theorem~\ref{Kachmar-LmV.8-V.9} and the hypothesis
$\lambda-\Theta(\widetilde \eta)<\zeta_0h^{2\rho}$, we obtain,
\begin{align*}
N\left
(h\lambda,\,\widetilde L_{h,\beta,S}^{\alpha,\eta}\right)
={\rm Card}\,\left(\{n\in\mathbb Z~;~
\mu_1\left(H_{h,\beta,2\pi n h^{1/2}S^{-1}}^{\alpha,\eta}\right)
\leq \lambda\}\right).
\end{align*}
Again, Theorem~\ref{Kachmar-LmV.8-V.9} yields the existence of a positive
constant $\zeta$ (that we may choose sufficiently large as we wish)
and a function $h\mapsto\widetilde\epsilon(h)$
such that, upon defining the subsets
\begin{align*}
\begin{split}
&\mathcal S_\pm=\bigg{\{}
n\in\mathbb Z~:~\left|2\pi nh^{1/2}S^{-1}
-\xi(\widetilde\eta)\right|\leq \zeta
h^{\rho},\\
&\quad\Theta(\widetilde\eta)+d_2(\alpha,\eta)\left(2\pi
nh^{1/2}S^{-1}-\xi(\widetilde\eta)\right)^2
-d_3(\alpha,\eta)\beta h^{1/2}\pm
h^{1/2}\widetilde\epsilon(h)<\lambda\bigg{\}},
\end{split}
\end{align*}
one gets the inclusion,
\begin{equation}\label{EquivalentCondition}
\mathcal S_+\subset
\bigg{\{}n\in\mathbb Z~;~
\mu_1\left(H_{h,\beta,2\pi n h^{1/2}S^{-1}}^{\alpha,\eta}\right)
\leq \lambda\bigg{\}}\subset \mathcal S_-.
\end{equation}
Therefore, we deduce that
$${\rm Card}\,\mathcal S_+
\leq
N\left
(h\lambda,\,\widetilde L_{h,\beta,S}^{\alpha,\eta}\right)\leq {\rm
  Card}\,S_-.
$$
On the other hand, thanks to (\ref{Hyp-lambda}), we may choose
$\zeta>0$ sufficiently large so that
\begin{eqnarray*}
&&\Theta(\widetilde\eta)+d_2(\alpha,\eta)\left(2\pi
nh^{1/2}S^{-1}-\xi(\widetilde\eta)\right)^2
-d_3(\alpha,\eta)\beta h^{1/2}\pm
h^{1/2}\widetilde\epsilon(h)<\lambda\\
&&\implies
\left|2\pi nh^{1/2}S^{-1}
-\xi(\widetilde\eta)\right|\leq \zeta
h^{\rho}\,.\end{eqnarray*}
With this choice of $\zeta$, one can rewrite $\mathcal S_\pm$ in the following
equivalent form
\begin{align*}
\begin{split}
&\mathcal S_\pm=\bigg{\{}
n\in\mathbb Z~:~ d_2(\alpha,\eta)\left(
2\pi nS^{-1}-h^{-1/2}\xi(\widetilde\eta)\right)^2\\
&\quad\quad\quad
\leq h^{-1/2}\left(d_3(\alpha,\eta)\beta
 +h^{-1/2}[\lambda-\Theta(\widetilde \eta)]
\pm
\widetilde\epsilon(h)\right)_+\bigg{\}},
\end{split}
\end{align*}
from which  one obtains a positive function $\epsilon_0(h)\ll1$ such
that
\begin{align*}
\begin{split}
&\left|
{\rm Card}\,\mathcal S_\pm
- \frac{h^{-1/4}S}{\pi\sqrt{d_2(\alpha,\eta)}}\sqrt{\left(
d_3(\alpha,\eta)\beta+h^{-1/2}[\lambda-\Theta(\widetilde\eta)]\right)_+}
\right|\\
&\leq \frac{h^{-1/4}S}{\pi\sqrt{d_2(\alpha,\eta)}}\sqrt{\left(
d_3(\alpha,\eta)\beta+h^{-1/2}[\lambda-\Theta(\widetilde\eta)]\right)_+}
\,\epsilon_0(h),
\end{split}
\end{align*}
when $S$ varies in a bounded interval $]-M,M[$.
This finishes the proof of the proposition.\hfill$\Box$\\

\subsection{The model operator on a Dirichlet strip}\label{CV-MO-D}
We continue to work  in the framework of the previous subsection by
keeping our choice of parameters $\beta$, $\eta$,
$\alpha\geq\frac12$, $h>0$, $\delta\in]\frac14,\frac12[$ and $S$.
Let us consider the operator $L_{h,\beta,S}^{\alpha,\eta}$ obtained
from $\widetilde L_{h,\beta,S}^{\alpha,\eta}$ by imposing additional
Dirichlet boundary conditions at $s\in\{0,S\}$, i.e.
$$L_{h,\beta,\eta}^{\alpha,\eta}~:~D(L_{h,\beta,S}^{\alpha,\eta})\ni
u\mapsto \widetilde L_{h,\beta,S}^{\alpha,\eta}u$$ with
$$D(L_{h,\beta,S}^{\alpha,\eta})=\{u\in D(\widetilde
L_{h,\beta,\eta}^{\alpha,\eta})~:~u(0,\cdot)=u(S,\cdot)=0\}.$$
Actually, $L_{h,\beta,S}^{\alpha,\eta}$ is the self adjoint operator
in $L^2\left(]0,S[\times]0,h^\delta[;\,(1-\beta h^\delta)\,\md s\md
t\right)$ associated with the quadratic form,
\begin{eqnarray*}
 Q_{h,\beta,S}^{\alpha,\eta}(u)&=&\int_0^S\int_0^{h^\delta}
\left(|h\partial_tu|^2+(1+2\beta t)\left|\left(h\partial_s+t
-\beta\frac{t^2}2\right)u\right|^2\right)\\
&&\hskip2cm\times(1-\beta t)\,\md s\md
t+h^{1+\alpha}\eta\int_0^S|u(s,0)|^2\,\md s,
\end{eqnarray*}
defined for functions $u$ in the form domain
$$D\left(Q_{h,\beta,S}^{\alpha,\eta}\right)=\{u\in H^1\left(
]0,S[\times]0,h^\delta[\right)~:~u(\cdot,h^\delta)=u(0,\cdot)
=u(S,\cdot)=0\}.$$  Using the same reasoning as that for the proof
of Lemma~\ref{JMP07-lem-N-MO'}, we get in the next lemma an estimate
of the spectral counting function of the operator
$L_{h,\beta,S}^{\alpha,\eta}$.
\begin{lem}\label{JMP-RF-prop4.5}
For each $M>0$, there exist constants $C>0$ and $h_0>0$ such that,
for all
$$\beta,\eta,S\in]-M,M[,\quad \varepsilon_0\in]0,S/2[,\quad \lambda\in\mathbb R,\quad
h\in]0,h_0[,$$ one has the estimate~:
$$\frac12 N\left(\lambda-C\varepsilon_0^{-2}h^2,\,\widetilde
L_{h,\beta,2(S-\varepsilon_0)}^{\alpha,\eta}\right)\leq
N\left(\lambda,\,L_{h,\beta,S}^{\alpha,\eta}\right)\leq
N\left(\lambda,\widetilde L_{h,\beta,S}^{\alpha,\eta}\right).$$
\end{lem}

\subsection{Spectral counting function in general
domains}\label{CE-GDom} We return in this subsection to the case of
a general smooth domain $\Omega$ whose boundary is compact.
\subsubsection*{An energy estimate}\ \\
Let us recall the notation that $\kappa_{\rm r}$ denotes the scalar
curvature of the boundary $\partial\Omega$. As was first noticed in
\cite{HeMo}, since the magnetic field is constant, the quadratic
form (\ref{VI-qform0}) can be estimated with a high precision by
showing the influence of the scalar curvature. This is actually the
content of the next lemma, which we quote from \cite[Lemma~4.7]{Fr}.
Before stating the estimate, let us recall that to a given function
$u\in H^1_{\rm loc}(\Omega)$, we associate by means of boundary
coordinates a function $\widetilde u$, see (\ref{VI-utilde}).\\

\begin{lem}\label{JMP-RF-Lemm4.7}
Let $\delta\in]\frac14,\frac12[$. There exists a constant $C>0$, and
for all
$$S\in[0,|\partial\Omega|[,\quad
\widetilde S\in]0,S[,$$ there exists a function $\phi\in
C^2([0,S]\times[0,Ch^\delta];\mathbb R)$ such that, for all $
\varepsilon\in[Ch^\delta,1]$, and for all $u\in H_{h,A}^1(\Omega)$
satisfying
$${\rm supp}\,\widetilde u\subset [0,S]\times[0,Ch^\delta],$$
one has the following estimate,
\begin{eqnarray*}
&&\left|q_{h,\Omega}(u)-Q_{h,\widetilde \kappa,S}\left(
e^{i\phi/h}\widetilde
u\right)\right|\\
&&\quad\leq C\left(h^\delta S\, Q_{h,\widetilde \kappa,S}\left(
e^{i\phi/h}\widetilde
u\right)+(h^{2+\delta}+Sh^{3\delta})\left\|e^{i\phi/h}\widetilde
u\right\|_{L^2}^2\right).
\end{eqnarray*}
Here $\widetilde\kappa=\kappa_{\rm r}(\widetilde S)$, and the
function $\widetilde u$ is associated to $u$ by
(\ref{VI-utilde}).\end{lem}

Let us mention that we omit $\alpha$ and $\eta$ from the notation in
Subsection~\ref{CV-MO-D} when $\eta\equiv0$.

\subsubsection*{Estimates of the counting function}\ \\
 As in
Section~\ref{proofThm1}, we introduce a partition of a thin
neighborhood of $\partial\Omega$~:\\
Given $N\in\mathbb N$ (that will be chosen later as a function of
$h$) such that
$$N=N(h)\gg 1\quad (h\to0),$$
we put
$$S=\frac{|\partial\Omega|}{N},\quad s_n=nS,\quad
\kappa_n=\kappa_{\rm r}(s_n),\quad n\in\{0,1,\cdots,N\}.$$ By this
way, we are able to estimate the spectral counting function of the
operator $P_{h,\Omega}^{\alpha,\gamma}$ by those of the operators
$L_{h,\kappa_n,S}^{\alpha,\gamma}$.

\begin{prop}\label{JMP-RF-prop4.8}
Let $\delta\in]\frac14,\frac12[$,
$\gamma_0=\displaystyle\min_{x\in\partial\Omega}\gamma(x)$ and
$\lambda=\lambda(h)$ such that
\begin{equation}\label{Hyp-lambda*}
\left|\lambda-\Theta\left(h^{\alpha-1/2}\gamma_0\right)\right|\ll1\quad
(h\to0).\end{equation} There exist constants $C>0$ and $h_0>0$ such
that, for all
$$h\in]0,h_0[,\quad \varepsilon_0\in]0,S/2[,\quad \widetilde
S_n\in[s_{n-1},s_n],\quad n\in\{1,\cdots,N\},$$ one has the estimate
\begin{align*}
\begin{split}
&\frac12\sum_{n=1}^N
N\left(h\lambda-C(Sh^{3\delta}+\varepsilon_0^{-2}h^2),\,\widetilde
L_{h,\widetilde\kappa_n,2(S+\varepsilon_0)}^{\alpha,\widetilde
\gamma_n}\right)\\
&\quad\leq
N\left(h\lambda,\,P_{h,\Omega}^{\alpha,\gamma}\right)\\
&\quad\leq \sum_{n=1}^N
N\left(h\lambda+C(Sh^{3\delta}+\varepsilon_0^{-2}h^2),\,\widetilde
L_{h,\widetilde\kappa_n,S+2\varepsilon_0}^{\alpha,\widehat\gamma_n}\right).
\end{split}
\end{align*}
Here
$$\widetilde\kappa_n=\kappa_{\rm r}(\widetilde S_n),\quad
\widetilde\gamma_n=\frac{\gamma(\widetilde S_n)+CS}{1+Ch^\delta
S},\quad \widehat\gamma_n=\frac{\gamma(\widetilde
S_n)-CS}{1-Ch^\delta S}.$$
\end{prop}
The proof is exactly as that of Proposition~\ref{VI-FrProp3.6} and
we omit it. Let us only mention the main points. Thanks to a
partition of unity associated with the intervals $[s_{n-1},s_n]$ and
the variational principle, the result follows from local estimates
of the quadratic form. For this sake, the procedure  consists of the
implementation of Lemma~\ref{JMP-RF-Lemm4.7}, bounding $\gamma(s)$
from above and below respectively by  $(1+Ch^\delta
S)\widetilde\gamma_n$ and $(1-Ch^\delta S)\widehat\gamma_n$  in each
$[s_{n-1},s_n]$ (thus getting  errors of order $S$) and finally of
the application of Lemma~\ref{JMP-RF-prop4.5}.

\subsubsection*{An asymptotic formula of the counting function}\ \\
Let us take  constants $\zeta_0>0$, $c_0>0$,
$\delta\in]\frac14,\frac12[$, and let us recall that we introduce a parameter
$\rho$ such that~:
$$0<\rho<\delta-\frac14\quad{\rm if~}\alpha=\frac12,\quad
0<\rho<\min\left(\delta-\frac14,\alpha-\frac12\right)\quad{\rm if~}\alpha>\frac12.$$
We take $\widetilde h_0>0$ and $\lambda=\lambda(h)$
such that,
\begin{equation}\label{Hyp-CE-conc}
c_0h^{1/2}
\leq
\left|\lambda-\Theta\left(h^{\alpha-1/2}\gamma_0\right)\right|<
\zeta_0h^{2\rho},\quad\forall~h\in]0,\widetilde
h_0].\end{equation}
Here
$$\gamma_0=\min_{x\in\partial\Omega}\gamma(x)$$ and $\gamma\in
C^\infty(\partial\Omega;\mathbb R)$ is the function in de\,Gennes'
boundary condition that we impose on functions in the domain of the
operator $P_{h,\Omega}^{\alpha,\gamma}$, see (\ref{JMP07DOp}).\\

\begin{thm}\label{JMP-CE-conc}
Let $\delta\in]\frac{5}{12},\frac12[$. With the above notations, we
have the following asymptotic formula as $h\to0$,
\begin{align*}
\begin{split}
&N\left(h\lambda,\,P_{h,\Omega}^{\alpha,\gamma}\right)
=\left(\int_{\partial\Omega}
\frac{h^{-1/4}}{\pi\sqrt{d_2(\alpha,\gamma(s))}}\,\times\right.\\
&\quad\left. \sqrt{\bigg{(} d_3(\alpha,\gamma(s))\kappa_{\rm
r}(s)+h^{-1/2}\left[\lambda-\Theta\left(h^{\alpha-1/2}\gamma(s)\right)\right]\bigg{)}_+}
\,\,\md s\right)\big{(}1+o(1)\big{)}.
\end{split}
\end{align*}
Here, for a given $(\alpha,\eta)\in\mathbb R\times\mathbb R$, the
parameters $d_2(\alpha,\eta)$ and $d_3(\alpha,\eta)$ have been
introduced in (\ref{CE-d2-d3}).
\end{thm}
\paragraph{\bf Proof.}
The proof is similar to that of the asymptotic formula
(\ref{JMP07-proofthm1}). There we have shown how to establish a
lower bound, so we show here how to establish an upper bound.\\
For each $n\in\{1,\cdots,N\}$, let us introduce the function (see
Lemma~\ref{JMP-RF-prop4.5})~:
$$f_n(\delta,S,\varepsilon_0)=d_3(\alpha,\widehat\gamma_n)
\widetilde\kappa_n+h^{-1/2}\left[\lambda+C(Sh^{3\delta-1}+\varepsilon_0^{-2}h)-\Theta\left(h^{\alpha-1/2}
\widehat\gamma_n\right)\right].$$ Here $\widehat\gamma_n$ and
$\widetilde\kappa_n$ are given by Proposition~\ref{JMP-RF-prop4.8}.
Then, combining Propositions~\ref{Fr-prop4.1} and
\ref{JMP-RF-prop4.8}, we get,
\begin{equation}\label{conc-proof1}
N\left(h\lambda,\,P_{h,\Omega}^{\alpha,\gamma}\right)\leq
\left(\sum_{n=1}^N\frac{h^{-1/4}(S+2\varepsilon_0)}{\pi\sqrt{d_2(\alpha,\widehat\gamma_n)}}\,\sqrt{[f_n(\delta,S,\varepsilon_0)]_+}\right)
\big{(}1+\epsilon_0(h)\big{)},
\end{equation}
where the function $\epsilon_0$ is independent of $N$ and satisfies
$$\lim_{h\to0}\epsilon_0(h)=0.$$
We recall also that $S=\frac{|\partial\Omega|}{N}$. We make the
following choice of $\varepsilon_0\in]0,S/2[$~:
$$\varepsilon_0=S^{1+\varsigma}\quad{\rm with}\quad\varsigma>0.$$
Then we pose the following condition on $S$ as $h\to0$,
$$Sh^{3\delta-1}+S^{-2-2\varsigma}h\ll
\left|\lambda-\Theta\left(h^{\alpha-1/2}
\widehat\gamma_n\right)\right|,\quad\forall~n\in\{1,\cdots,N\}.$$
By  the hypothesis in (\ref{Hyp-CE-conc}), it suffices to choose $S$
in the following way~:
$$Sh^{3\delta-1}+S^{-2-2\varsigma}h\ll h^{1/2}\quad (h\to0).$$
This yields,
$$h^{\frac1{4(1+\varsigma)}}\ll S\ll h^{3(\frac12-\delta)}$$
and we notice that a choice of $\varsigma>0$ such that
$$\frac1{4(1+\varsigma)}>3\left(\frac12-\delta\right)$$
is possible only if $\delta\in]\frac5{12},\frac12[$ (this will yield
when $\alpha=\frac12$ that $\rho_0\in]\frac16,\frac14[$, $\rho_0$
    being introduced in Theorem~\ref{Kachmar-LmV.8-V.9}). With this
choice, the upper bound (\ref{conc-proof1}) becomes (for a possibly
different $\widetilde\epsilon_0(h)\ll1$),
\begin{equation}\label{conc-proof2}
N\left(h\lambda,\,P_{h,\Omega}^{\alpha,\gamma}\right)\leq
\big{(}1+\widetilde\epsilon_0(h)\big{)}
\sum_{n=1}^N\frac{h^{-1/4}S}{\pi\sqrt{d_2(\alpha,\widehat\gamma_n)}}\,
\sqrt{[g_n(\lambda,\alpha)]_+} ,
\end{equation}
with
$$g_n(\lambda,\alpha)=d_3(\alpha,\widehat\gamma_n)
\widetilde\kappa_n+h^{-1/2}\left[\lambda-\Theta\left(h^{\alpha-1/2}
\widehat\gamma_n\right)\right].$$
Replacing $\widehat\gamma_n$ by $\gamma_n=\gamma(\widetilde S_n)$ in
(\ref{conc-proof2}) will yield an error of the order $\mathcal O(S)$,
and
the sum on the right hand side of
(\ref{conc-proof2}) becomes a Riemann sum. We therefore conclude
the following upper
bound $$N\left(h\lambda,\,P_{h,\Omega}^{\alpha,\gamma}\right)\leq
\big{(}1+\widetilde\epsilon_0(h)\big{)}
\int_{\partial\Omega}\frac{h^{-1/4}}{\pi\sqrt{d_2(\alpha,\gamma(s))}}\,
\sqrt{[g(\lambda,\alpha;s)]_+}\,\md s,$$ with
$$g(\lambda,\alpha;s)=d_3(\alpha,\gamma(s))
\widetilde\kappa(s)+h^{-1/2}\left[\lambda-\Theta\left(h^{\alpha-1/2}
\gamma(s)\right)\right].$$ By a similar argument, we get a lower
bound. \hfill$\Box$\\

\begin{rem}\label{relax-hyp}
{\rm When relaxing the hypotheses of Theorem~\ref{JMP-CE-conc} by
allowing $\lambda$ to satisfy (compare with (\ref{Hyp-CE-conc})):
$$\left|\lambda-\Theta\left(h^{\alpha-1/2}\gamma_0\right)\right|=o(h^{1/2})
\quad{\rm as}~h\to0\,,$$
the result for the counting function becomes (as can be checked by
adjusting the proof of Theorem~\ref{JMP-CE-conc})}
\begin{align*}
\begin{split}
&N\left(h\lambda,\,P_{h,\Omega}^{\alpha,\gamma}\right)
=\left(\int_{\partial\Omega}
\frac{h^{-1/4}}{\pi\sqrt{d_2(\alpha,\gamma(s))}}\,
\sqrt{\bigg{(} d_3(\alpha,\gamma(s))\kappa_{\rm
r}(s)\bigg{)}_+}
\,\,\md s\right)\big{(}1+o(1)\big{)}.
\end{split}
\end{align*}
\end{rem}

\text{ }\vskip0.3cm
\subsubsection*{Proof of Theorem~\ref{JMP07thm2}}\ \\
We recall in this case that $\frac12<\alpha<1$ and that
$$\lambda(h)=\Theta_0+3aC_1h^{\alpha-\frac12}\quad{\rm with}\quad
a\in\mathbb R\setminus\{\gamma_0\}.$$
Here $C_1>0$ is the universal constant introduced in
(\ref{JMP07C1}).\\
In this specific regime, (\ref{Hyp-CE-conc}) is verified when making a
choice of\break
$\rho\in]0,\frac12\min(\delta-\frac14,\alpha-\frac12)[$\,.\\
The leading order term of the
integrand in the
asymptotic formula of Theorem~\ref{JMP-CE-conc} is, up to a
multiplication by a positive constant,
$$\sqrt{h^{-1/2}\left[\lambda-\Theta\left(h^{\alpha-1/2}\gamma(s)\right)\right]_+}~.$$
We write by using the asymptotic expansion of $\Theta(\cdot)$ given by Taylor's
formula (see (\ref{Th'-gam})-(\ref{C1}))~:
$$\Theta\left(h^{\alpha-1/2}\gamma(s)\right)=\Theta_0+3C_1\gamma(s)h^{\alpha-1/2}+\mathcal
O(h^{2\alpha-1}),\quad (h\to0).$$ Therefore, it results from
Theorem~\ref{JMP-CE-conc},
$$N\left(h\lambda,\,P_{h,\Omega}^{\alpha,\gamma}\right)
=\left(\frac{h^{\frac12(\alpha-\frac32)}}{\pi\sqrt{\xi_0}}\int_{\partial\Omega}\sqrt{
\big{(}a-\gamma(s)\big{)}_+}\,\md s\right)(1+o(1)).
$$
When $a=\gamma_0$, we may encounter the regime of
Remark~\ref{relax-hyp}, hence by using the result of that remark and
noticing that when $\frac12<\alpha<1$
$$h^{-1/4}\ll h^{\frac12(\alpha-\frac32)}\ll h^{-1/2}\quad{\rm
  as~}h\to0\,,$$
we recover the asymptotic expansion announced in
  Theorem~\ref{JMP07thm3} in the present case.
\hfill$\Box$

\subsubsection*{Proof of Theorem~\ref{JMP07thm3}}\ \\
In this case $\alpha=\frac12$ and
$$h^{1/2}\ll|\lambda-\Theta(\gamma_0)|\leq\zeta_0h^{\varrho}\quad
{\rm with}\quad 0<\varrho<\frac12.$$
Taking $\rho=\varrho/2$, then we may choose
$\delta\in]\frac5{12},\frac12[$  such that (\ref{Hyp-CE-conc}) is
    satisfied. Thus,
 the asymptotic formula of
Theorem~\ref{JMP-CE-conc} is still valid in this regime, and the
leading order term of the integrand is, up to a multiplicative
constant,
$$\frac{h^{-1/2}}\pi
\sqrt{\frac{\left[\lambda-\Theta\left(\gamma(s)\right)\right]_+}{d_2\left(\frac12,\gamma(s)\right)}}\,.$$
This proves the theorem.\hfill$\Box$

\subsubsection*{Proof of Theorem~\ref{JMP07thm4}}\ \\
Again, the proof follows from Theorem~\ref{JMP-CE-conc} and the
properties of the function $\Theta(\cdot)$.\hfill$\Box$\\

\section*{Acknowledgements} The author would like to express his thanks
to R. Frank for the fruitful discussions around the subject, and
also to B.~Helffer for his helpful remarks.
He wishes also to thank
the anonymous referees for their careful reading of the paper, for
pointing out many corrections and for their many helpful
suggestions.

\appendix
\section{Boundary coordinates}

We recall now
the definition of the standard coordinates that straightens a
portion of the boundary $\partial\Omega$. Given $t_0>0$, let us
introduce the following neighborhood of the boundary,
\begin{equation}\label{App:N-t0}
\mathcal N_{t_0}=\{x\in\mathbb R^2;\quad {\rm
dist}(x,\partial\Omega)<t_0\}.
\end{equation}
As the boundary is smooth, let
$s\in]-\frac{|\partial\Omega|}2,\frac{|\partial\Omega|}2]\mapsto
M(s)\in\partial\Omega$ be a regular parametrization of
$\partial\Omega$ that satisfies~:
$$\left\{
\begin{array}{l}
s \text{ is
the oriented `arc length' between }M(0) \text{ and }
M(s).\\
T(s):=M'(s) \text{ is a unit tangent vector
to }\partial\Omega\text{ at the point }M(s).\\
\text{The orientation is positive, i.e. }{\rm det}(T(s),\nu(s))=1.
\end{array}
\right. $$ We recall that $\nu(s)$ is the unit outward normal of
$\partial\Omega$ at the point $M(s)$. The scalar curvature
$\kappa_{\rm r}$ is now defined by~:
\begin{equation}\label{App:kappa-r}
T'(s)=\kappa_{\rm r}(s)\nu(s).
\end{equation}
When $t_0$ is sufficiently small, the map~:
\begin{equation}\label{App:Phi(s,t)}
\Phi:\,]-|\partial\Omega|/2,|\partial\Omega|/2]\times
]-t_0,t_0[\,\ni(s,t)\mapsto M(s)-t\nu(s)\in \mathcal N_{t_0},
\end{equation}
is a diffeomorphism. For $x\in\mathcal N_{t_0}$, we write,
\begin{equation}\label{App:Phi-1}
\Phi^{-1}(x):=(s(x),t(x)),
\end{equation}
where
$$t(x)={\rm dist}(x,\partial\Omega)\text{ if }x\in\Omega \quad
\text{and }t(x)=-{\rm dist}(x,\partial\Omega)\text{ if
}x\not\in\Omega.$$ The Jacobin of the transformation $\Phi^{-1}$ is
equal to,
\begin{equation}\label{App:Jac}
a(s,t)={\rm det}\left(D\Phi^{-1}\right)=1-t\kappa_{\rm r}(s).
\end{equation}
To a vector field $A=(A_1,A_2)\in H^1(\mathbb R^2;\mathbb R^2)$, we
associate the vector field
$$\tilde A=(\tilde A_1,\tilde A_2)\in
H^1(]-|\partial\Omega|/2,|\partial\Omega|/2]\times]-t_0,t_0[;\mathbb
R^2)$$
by the following relations~:
\begin{equation}\label{App:chmnew}
\tilde A_1(s,t)=(1-t\kappa_{\rm r}(s)) \vec A(\Phi(s,t))\cdot
M'(s),\quad \tilde A_2(s,t)=\vec A(\Phi(s,t))\cdot\nu(s).
\end{equation}
We get then the following change of variable formulae.
\begin{proposition}\label{App:transf}
Let $u\in H^1_A(\mathbb R^2)$ be supported in $\mathcal N_{t_0}$.
Writing $\widetilde u(s,t)=u(\Phi(s,t))$, then we have~:
\begin{equation}\label{App:qfstco}
\int_{\Omega}\left|(\nabla-iA)u\right|^2dx=
\int_{-\frac{|\partial\Omega|}2}^{\frac{|\partial\Omega|}2}
\int_{0}^{t_0}\left[ |(\partial_s-i\tilde A_1)\widetilde
u|^2+a^{-2}|(\partial_t-i\tilde A_2)\widetilde u|^2\right]a\,dsdt,
\end{equation}
\begin{equation}\label{App:qfstco'}
\int_{\Omega^c}\left|(\nabla-iA)u\right|^2dx=
\int_{-\frac{|\partial\Omega|}2}^{\frac{|\partial\Omega|}2}
\int_{-t_0}^{0}\left[ |(\partial_s-i\tilde A_1)\widetilde
u|^2+a^{-2}|(\partial_t-i\tilde A_2)\widetilde u|^2\right]a\,dsdt,
\end{equation}
and
\begin{equation}\label{App:nostco}
\int_{\mathbb R^2}
|u(x)|^2\,dx=\int_{-\frac{|\partial\Omega|}2}^{\frac{|\partial\Omega|}2}
\int_{-t_0}^{t_0} |\widetilde u(s,t)|^2a\,dsdt.
\end{equation}
\end{proposition}

We have also the relation~:
$$\left(\partial_{x_1}A_2-\partial_{x_2}A_1\right)dx_1\wedge dx_2=
\left(\partial_s\tilde A_2-\partial_t\tilde A_1\right)a^{-1}ds\wedge
dt,$$ which gives,
$${\rm curl}_{(x_1,x_2)}\, A=\left(1-t\kappa_{\rm r}(s)\right)^{-1}{\rm curl}_{(s,t)}\,\tilde A.$$

We give in the next proposition a standard choice of gauge.

\begin{prop}\label{App:Agd1}
Consider a vector field $A=(A_1,A_2)\in C^1_{\rm loc}(\mathbb
R^2;\mathbb R^2)$ such that
$${\rm curl }\,A=1\quad \text{ in }\mathbb R^2.$$ For each point $x_0\in\partial\Omega$, there
exist a neighborhood $\mathcal V_{x_0}\subset\mathcal N_{t_0}$ of
$x_0$ and a smooth real-valued function $\phi_{x_0}$ such that the
vector field $A_{new}:=A-\nabla\phi_{x_0}$ satisfies in $\mathcal
V_{x_0}$~:
\begin{equation}\label{App:HeMo'}
\tilde A_{new}^2=0,
\end{equation}
and,
\begin{equation}\label{App:HeMo}
\tilde A_{new}^1=-t\left(1-\frac t2\kappa_{\rm r}(s)\right),
\end{equation}
with $\tilde A_{new}=(\tilde A_{new}^1,\tilde A_{new}^2)$.
\end{prop}

\end{document}